\documentclass[preprint,11pt]{elsarticle}
\usepackage[margin=1in]{geometry}
\usepackage{amsmath,amssymb,amsthm}
\usepackage{enumitem}
\usepackage{tikz}
\usepackage{booktabs}
\usepackage{hyperref}
\usetikzlibrary{arrows.meta,positioning,bending}
\hypersetup{colorlinks=true,linkcolor=blue!50!black,citecolor=blue!50!black,urlcolor=blue!50!black}

\numberwithin{equation}{section}
\newtheorem{theorem}[equation]{Theorem}
\newtheorem{proposition}[equation]{Proposition}
\newtheorem{lemma}[equation]{Lemma}

\newtheorem{definition}[equation]{Definition}
\newtheorem{remark}[equation]{Remark}

\newcommand{\Aut}{\operatorname{Aut}}

\newcommand{\MM}{\mathcal{M}}

\newcommand{\ZZ}{\mathbb{Z}}
\newcommand{\RR}{\mathbb{R}}
\newcommand{\CC}{\mathbb{C}}
\newcommand{\PP}{\mathbb{P}}
\newcommand{\TT}{\mathcal{T}}

\newcommand{\WM}{Whitehead move}
\newcommand{\interior}{\operatorname{int}}
\newcommand{\rank}{\operatorname{rank}}

\makeatletter
\def\ps@pprintTitle{%
  \let\@oddhead\@empty
  \let\@evenhead\@empty
  \def\@oddfoot{\reset@font\hfil\thepage\hfil}%
  \let\@evenfoot\@oddfoot}
\makeatother

\begin{document}

\begin{frontmatter}

\title{Delaunay Combinatorics of the Moduli and Teichm\"uller Space
of Closed Genus-Two Riemann Surfaces}

\author{Armando Rodado Amaris}
\ead{armando@rodadolab.com}
\affiliation{organization={Centre for Quantum Software and Information, University of Technology Sydney},
  city={Sydney},
  state={NSW},
  country={Australia}}

\begin{abstract}
Every closed Riemann surface of genus two is hyperelliptic, so it is
determined by six branch points on the sphere, which become cone points
of angle $\pi$ in the quotient orbifold metric. The Delaunay
decomposition of these cone points, dual to the Voronoi diagram of the
Weierstrass points, gives a trivalent graph with eight vertices, twelve
edges and six faces. We classify these graphs. When the dual
triangulation is simple there are exactly two types, dual to the
octahedron and to a stacked triangulation, with orientation-preserving
automorphism groups of orders $24$ and $2$, and we determine the
Whitehead-move (diagonal-flip) orbit connecting them. Allowing bigons,
the orbit closes at exactly twenty types. Iterated edge contraction from
the two simple types gives a seven-level hierarchy of $2,6,22,32,22,5$
and $1$ types, ending at two vertices joined by six parallel edges;
closing under both contraction and expansion gives all $280$ maps on the
sphere with six vertices and no vertex of degree one. The rank of the
linear cone-angle system of each type equals $\dim_\RR\TT_2=6$. We
realize both simple types by explicit hyperbolic cone spheres, in closed
form for the octahedral type, build the corresponding genus-two surfaces
by cutting and doubling, and show that the boundary $x_e=0$ of the angle
system is exactly the concyclicity of four branch points. The
classification, the automorphism groups and the flip orbit are formally
verified in Lean~4, and every enumeration is reproduced by an
independent program, which also determines which types are realized by
Delaunay decompositions.
\end{abstract}


\end{frontmatter}

\tableofcontents

\section{Introduction}

\subsection{Background}
\label{sec:background}

The moduli space $\MM_g$ and Teichm\"uller space $\TT_g$ of closed
Riemann surfaces of genus $g\ge 2$ have been studied combinatorially
by several essentially different routes. Harer's cell decomposition
of the decorated moduli space via Jenkins--Strebel differentials and
ribbon graphs~\cite{Harer1986,Kontsevich1992} realizes $\MM_{g,n}$ as
a union of cells indexed by trivalent ribbon graphs. Independently,
Penner's decorated Teichm\"uller theory~\cite{Penner1987} and the
Bowditch--Epstein construction~\cite{BowditchEpstein1988} attach to
every decorated hyperbolic surface with cusps a canonical
Delaunay/Voronoi cell decomposition, giving $\TT_{g,n}\times
\RR_{>0}^n$ the structure of a cell complex indexed by ideal
triangulations, with adjacent cells related by diagonal flips.
Rivin~\cite{Rivin1994} and Thurston~\cite{ThurstonShapes1998},
building on Andreev's theorem, developed the closely related theory
of circle patterns realizing a prescribed triangulation with
prescribed intersection angles, which is the Euclidean/cone-metric
analogue of the same idea.

The present paper specializes this circle of ideas to genus two,
where an extra piece of rigid structure is available: every closed
Riemann surface of genus two is hyperelliptic, with a \emph{unique}
hyperelliptic involution~\cite{FarkasKra1992}, so that its
Teichm\"uller space is canonically identified with the Teichm\"uller
space of the six-times-marked sphere. This lets the
Bowditch--Epstein/Rivin--Thurston circle-pattern picture be pushed
down from the surface to the sphere, where the combinatorics becomes
completely explicit and, as we show, completely enumerable: the
possible combinatorial types of the resulting trivalent graph, and
the \WM{} moves connecting them, form a finite and fully classified
structure. This is analogous in spirit to the finite classification
of ideal triangulations underlying the (once-punctured torus)
case classically studied via the Farey tessellation, but genus two is
the first closed case where the construction is available, precisely
because of hyperellipticity.

This approach --- the graph associated with a genus-two surface
through the Voronoi diagram of its Weierstrass points, the
classification of these graphs by Whitehead moves, and the Delaunay
realization problem (DRP) of deciding which combinatorial types occur,
solved there with Springborn's circle-pattern existence
theorem~\cite{Springborn2003} --- originates in the author's 2007 PhD
thesis~\cite{Rodado2007}. The present paper revisits it, extends it,
and recomputes its combinatorial results by independent methods. Our
terminology differs from the thesis in places: our ``CE-graph'' is the
thesis's ``generic graph'', and our generic types $d_1,d_2$, singled
out by having a \emph{simple} dual triangulation, are the thesis's Gen1
and Gen2.

\subsection{Summary of results}

Section~\ref{sec:prelim} recalls the hyperelliptic quotient
construction and the combinatorial (dart-system) formalism used
throughout. Section~\ref{sec:classification} proves the main
classification theorem: a \emph{generic} CE-graph (one whose dual
triangulation is simple) is, up to embedded graph equivalence, one of
exactly two types, which we call the cube type $d_1$ and the stacked
type $d_2$ (Theorem~\ref{thm:classification}), with orientation-preserving
automorphism groups of order $24\cong S_4$ and $2$
(Proposition~\ref{prop:aut}). Section~\ref{sec:orbit} computes the
\WM{} orbit of $d_1$ completely (Theorem~\ref{thm:orbit}): every one
of its twelve facets is admissible and lands on $d_2$, while of
$d_2$'s twelve facets exactly one reverses to $d_1$, exactly four are
self-transitions forming a single orbit under the automorphism group of $T_2$, and the
remaining seven are inadmissible. Section~\ref{sec:teich} identifies
the target dimension six by two independent arguments and
Section~\ref{sec:circlepacking} shows that the rank of the natural
linear cone-angle system matches it exactly, via an elementary
argument (Lemma~\ref{lem:incidencerank}) valid uniformly for both
types. Section~\ref{sec:realizability} carries out the remaining
nonlinear step directly from the hyperbolic law of cosines, with no
appeal to any named circle-pattern existence theorem: it constructs
an explicit hyperbolic sphere realizing each generic type's Delaunay
decomposition (Propositions~\ref{prop:cubeconstruction}--\ref{prop:stackedconstruction}),
verified independently via a Gauss--Bonnet cross-check and a local
Delaunay criterion, builds the corresponding genus-two surface by
cutting and doubling (Theorem~\ref{thm:realizability}), and extends
this to one admissible edge contraction of each generic type
(Theorem~\ref{thm:contractrealize}), identifying the correct
geometric meaning of the linear system's boundary degeneracy along
the way (Proposition~\ref{prop:xezero}). Section~\ref{sec:beyond}
extends the classification in two
independent directions: Section~\ref{sec:fullorbit} drops genericity
within the trivalent level and shows (Theorem~\ref{thm:fullorbit})
that the full \WM{} orbit, allowing bigons, is exactly twenty types,
with face-degree signature shown not to be a complete invariant of
type; Section~\ref{sec:hierarchy} leaves the trivalent level
altogether via iterated edge contraction, revealing a complete
seven-level hierarchy of exactly $2,6,22,32,22,5,1$ combinatorial
types (Theorem~\ref{thm:hierarchy}) descending to a single minimal
graph, and identifies this hierarchy with the codimension
stratification of the Delaunay decomposition by concyclicity
degeneracies. Section~\ref{sec:atlas} displays all ninety types
explicitly, level by level, and Section~\ref{sec:regeneration} shows
the hierarchy can be regenerated upward from the minimal graph alone
by mixing expansion and contraction moves (Theorem~\ref{thm:regeneration}),
but not by expansion moves alone, which reach only $128$ of the
resulting $280$ types and miss twenty-four of the ninety generic-rooted
types outright (Proposition~\ref{prop:noexpansiononly}). The finite
combinatorics of Sections~\ref{sec:classification} and~\ref{sec:orbit}
is formally verified in Lean~4 (Remark~\ref{rem:lean}), and every
enumeration in the paper is reproduced by an independent program
(Remark~\ref{rem:indep}). Finally,
Section~\ref{sec:discussion} discusses the relation of these results
to the literature cited above and to the classical theory of
automorphisms of genus-two curves~\cite{Bolza1887,Broughton1991}, and
states what remains open toward extending realizability to the full
twenty-type orbit and the deeper contraction hierarchy.

\section{Preliminaries}
\label{sec:prelim}

\subsection{The hyperelliptic quotient}

\begin{lemma}[Every genus-two surface is hyperelliptic]
\label{lem:allhyp}
Let $X$ be a closed Riemann surface of genus $2$. The canonical map
$\phi_K:X\to\PP^1$ is a degree-two branched cover, and its deck
involution $\iota$ is the unique hyperelliptic involution of $X$.
\end{lemma}

\begin{proof}
The canonical bundle $K_X$ has degree $2g-2=2$, and by Riemann--Roch
$\ell(K_X)=g=2$, so $\phi_K$ maps $X$ to
$\PP^{\ell(K_X)-1}=\PP^1$. The canonical system is base-point-free
for every curve of genus $\ge 2$, so $\phi_K$ is a morphism of degree
$\deg K_X=2$; it cannot be an embedding, since an embedded curve in
$\PP^1$ would itself be $\PP^1$, of the wrong genus. Hence $\phi_K$
exhibits $X$ as a hyperelliptic double cover of $\PP^1$. Conversely, a
degree-two map $X\to\PP^1$ is given by a base-point-free pencil $|D|$
with $\deg D=2$ and $\ell(D)=2$. Riemann--Roch gives
$\ell(K_X-D)=\ell(D)-\deg D-1+g=1$, so $K_X-D$ is an effective divisor of
degree $0$ and $D\sim K_X$. Every degree-two map to $\PP^1$ is therefore
$\phi_K$ followed by an automorphism of $\PP^1$, and the hyperelliptic
involution is unique~\cite{FarkasKra1992}.
\end{proof}

By Riemann--Hurwitz, $2(2)-2 = 2(2\cdot 0 -2) + \#\mathrm{branch}$
forces exactly six branch points $p_1,\dots,p_6\in\PP^1$, the images
of the six fixed points of $\iota$. Writing $S=X/\langle\iota\rangle
\cong\PP^1$, the map $X\to S$ is unramified away from these six
points.

\begin{proposition}[Cone angle $\pi$]
\label{prop:coneangle}
Let $h$ be the hyperbolic metric on $X$ in its conformal class. Then
$\iota^*h=h$, and the induced metric on the orbifold $S=X/\langle
\iota\rangle$ has a cone point of angle $\pi$ at each $p_i$ and is
smooth elsewhere.
\end{proposition}

\begin{proof}
$\iota$ is a biholomorphism, so $\iota^*h$ is a hyperbolic metric in
the same conformal class as $h$; by uniqueness of the uniformizing
metric, $\iota^*h=h$. Near a fixed point of $\iota$, choose a
holomorphic coordinate $z$ in which the order-two automorphism acts
linearly, hence as $z\mapsto -z$; the quotient coordinate is
$w=z^2$, under which a full turn $\theta=2\pi$ upstairs covers only
angle $\pi$ downstairs.
\end{proof}

Proposition~\ref{prop:coneangle} says precisely that $S$, topologically
a sphere, carries an orbifold structure with six cone points of order
two (angle $2\pi/2=\pi$) and that this orbifold is hyperbolic. We make
this explicit, since it is the geometric fact the rest of the paper's
lifting arguments rest on.

\begin{proposition}[$S$ is a hyperbolic orbifold]
\label{prop:orbifold}
The orbifold $S$, with underlying topological space $S^2$ and a cone
point of order $2$ at each $p_i$, has negative orbifold Euler
characteristic,
\[
\chi^{\mathrm{orb}}(S) = \chi(S^2) - \sum_{i=1}^{6}\Big(1-\frac{1}{2}\Big)
= 2 - 3 = -1 < 0,
\]
so it admits a complete hyperbolic structure in its conformal class,
unique up to isometry by the orbifold uniformization
theorem~\cite{ThurstonNotes}; this is exactly the structure of
Proposition~\ref{prop:coneangle}. Its orbifold fundamental group is a
Fuchsian group
\[
\Gamma := \pi_1^{\mathrm{orb}}(S)
= \big\langle x_1,\ldots,x_6 \bigm| x_i^2=1,\ x_1x_2\cdots x_6=1
\big\rangle \ \le\ \mathrm{PSL}(2,\RR),
\]
where $x_i$ represents a loop once around $p_i$ and acts on the
universal cover $\mathbb H^2$ as the elliptic rotation by $\pi$ about
a lift $\tilde p_i$. The map $\Gamma\to\ZZ/2$, $x_i\mapsto1$, is a
well-defined homomorphism (the relation maps to $6\equiv0$), and its
kernel $\Gamma'\trianglelefteq\Gamma$ is torsion-free of index $2$
with $\mathbb H^2/\Gamma'\cong X$: the double cover $X\to S$
realizing $\iota$ is exactly the orbifold cover corresponding to
$\Gamma'\le\Gamma$, with $\iota$ the deck transformation of the
nontrivial coset $\Gamma\smallsetminus\Gamma'$.
\end{proposition}

\subsection{CE-graphs and dart systems}

\begin{definition}[CE-graph]
A \emph{CE-graph} is a graph $G$, connected and cellularly embedded on
$S^2$, such that every vertex has valence three, $G$ has six faces,
and no face is a monogon.
\end{definition}

The monogon exclusion is not an arbitrary normalization: for the
graph this paper's construction actually produces, it is a theorem
rather than a hypothesis, proved once the necessary orbifold
machinery is in place (Theorem~\ref{thm:nomonogon}).

Euler's formula with $3v=2e$ and $f=6$ forces $v=8$, $e=12$ for every
CE-graph. Loops, parallel edges, and contractions are handled
unambiguously via \emph{dart systems}: a dart system is a triple
$(D,\alpha,\rho)$ where $D$ is a finite set of darts, $\alpha:D\to D$
is a fixed-point-free involution pairing the two darts of each edge,
and $\rho:D\to D$ records the cyclic order of darts at each vertex
(so the cycles of $\rho$ are the vertices); the face permutation is
$\phi=\rho\circ\alpha$, whose cycles are the oriented face boundary
walks. Two embedded graphs are \emph{equivalent} if there is an
orientation-preserving homeomorphism of $S^2$ carrying one embedding
to the other; equivalently, a bijection $\varphi:D\to D'$ with
$\varphi\alpha=\alpha'\varphi$ and $\varphi\rho=\rho'\varphi$. This
notion is strictly finer than abstract graph isomorphism, which
discards $\rho$ and can identify combinatorially distinct
embeddings; all classification statements below are with respect to
embedded graph equivalence.

A CE-graph $G$ arises as the graph dual to a Delaunay decomposition
of the six cone points $p_1,\dots,p_6$ with respect to the metric of
Proposition~\ref{prop:coneangle}, pulled back to $S$ from the
hyperbolic metric on $X$ via the descent argument below
(Proposition~\ref{prop:descent}). This decomposition is a triangulation,
and $G$ is trivalent, exactly when no empty circle in the universal
cover passes through four lifts of the cone points
(Proposition~\ref{prop:descent}(c)). We call $G$ \emph{generic} if
moreover its dual triangulation is simple, that is, has no loops and no
multiple edges. A trivalent $G$ whose dual has a loop or a multiple edge
is not generic, even though no four points are concyclic; such graphs
occur for genuine surfaces (Section~\ref{sec:fullorbit}).

\section{Classification of generic CE-graphs}
\label{sec:classification}

\begin{lemma}[Voronoi--Delaunay duality]
\label{lem:duality}
Let $G$ be a generic CE-graph and $T=G^*$ its dual. Then $T$ has six
vertices, twelve edges, eight triangular faces, and is simple (no
loops, no multiple edges); consequently $T$ is a simplicial
triangulation of $S^2$ on the six marked points.
\end{lemma}

\begin{proof}
Vertex, edge and face counts dualize directly, trivalence of $G$
forces every face of $T$ to be a triangle, and $T$ is simple by the
definition of genericity.
\end{proof}

The base CE-graph definition excludes only monogons; it does not by
itself exclude a bigon (two faces sharing two edges, dual to a
multiple edge in $T$). Such a configuration is still a valid
CE-graph, simply not generic; Section~\ref{sec:fullorbit} shows that
dropping genericity enlarges the classification substantially, so
genericity is a real hypothesis and not a bookkeeping convenience.

\begin{lemma}[Degree bound and balance]
\label{lem:degrees}
For a simplicial triangulation $T$ of $S^2$ with exactly six
vertices, every vertex has degree $3$, $4$, or $5$, and if $x_k$ is
the number of degree-$k$ vertices, $x_3+x_4+x_5=6$ and $x_5=x_3$.
\end{lemma}

\begin{proof}
Simplicial triangulations with $\ge 4$ vertices have minimum degree
$\ge 3$; simplicity on six vertices forces maximum degree $\le 5$.
Euler's formula gives $e=12$, so $\sum \deg = 24 = \sum_k k\,x_k$;
subtracting $4\sum_k x_k = 24$ gives $x_5-x_3=0$.
\end{proof}

\begin{lemma}[Pentagon lemma]
\label{lem:pentagon}
Every triangulation of a pentagon by non-crossing diagonals, using no
vertices beyond the five boundary vertices, is a fan from a single
boundary vertex.
\end{lemma}

\begin{proof}
Two diagonals of a pentagon cross exactly when they do not share an
endpoint (checked directly on all $\binom{5}{2}=10$ pairs); a
triangulation uses exactly two non-crossing diagonals, which must
therefore share a common endpoint, and a pentagon has only two
diagonals through any given vertex.
\end{proof}

\begin{proposition}[The all-degree-four case]
\label{prop:octa}
There is exactly one six-vertex simplicial triangulation of $S^2$
with all degrees equal to four, up to isomorphism: the octahedral
triangulation, with $1$-skeleton $K_{2,2,2}$.
\end{proposition}

\begin{proof}
A $4$-regular simple graph on six vertices is the complement of a
perfect matching; all perfect matchings on six vertices are
isomorphic, so their complements are too, giving the unique class
$K_6\setminus(\text{perfect matching})=K_{2,2,2}$, the octahedron
graph.
\end{proof}

\begin{proposition}[The case with a degree-five vertex]
\label{prop:stack}
There is exactly one six-vertex simplicial triangulation of $S^2$
with a vertex of degree five, up to isomorphism, with degree sequence
$(3,3,4,4,5,5)$.
\end{proposition}

\begin{proof}
If $A$ has degree $5$ it is adjacent to all others, so its link is a
pentagon on the remaining five vertices, and the complementary disk
is a triangulated pentagon with no further vertices available; by
Lemma~\ref{lem:pentagon} this is a fan from some vertex $w$. Writing
the pentagon $u,w,x,y,z$ in cyclic order, the fan adds diagonals
$wy,wz$, giving degree $5$ at $w$, degree $4$ at $y,z$, and degree $3$
at $u,x$ --- forced sequence $(3,3,4,4,5,5)$ --- and uniqueness
follows since all five choices of fan vertex are related by the
pentagon's rotational symmetry.
\end{proof}

\begin{theorem}[Classification of generic CE-graphs]
\label{thm:classification}
Up to embedded graph equivalence there are exactly two generic
CE-graphs: the \emph{cube type} $d_1$, dual to the octahedral
triangulation, all of whose six faces are quadrilaterals; and the
\emph{stacked type} $d_2$, dual to the triangulation of
Proposition~\ref{prop:stack}, whose six faces have sizes
$3,3,4,4,5,5$.
\end{theorem}

\begin{proof}
By Lemma~\ref{lem:degrees}, a six-vertex triangulation either has no
degree-$5$ vertex (forcing the all-degree-four case,
Proposition~\ref{prop:octa}) or has one (Proposition~\ref{prop:stack},
whose proof also excludes $x_3\in\{1,3\}$); these cases are
exhaustive, mutually exclusive, and dualize to $d_1,d_2$ by
Lemma~\ref{lem:duality}.
\end{proof}

Both triangulations are $3$-connected --- a standard fact for
simplicial $2$-sphere triangulations with $\ge 4$ vertices --- so by
Steinitz's theorem~\cite{Grunbaum2003} they are boundary complexes of
convex simplicial $3$-polytopes, and by Whitney's
theorem~\cite{Whitney1932} each has a unique embedding in $S^2$ up to
reflection; consequently classifying generic CE-graphs up to embedded
equivalence is exactly classifying these triangulations up to
orientation-preserving isomorphism, with no information lost in
passing to the dual.

\begin{proposition}[Orientation-preserving automorphisms]
\label{prop:aut}
$\Aut^+(d_1)$ is the rotation group of the cube, of order $24\cong
S_4$. $\Aut^+(d_2)$ has order $2$.
\end{proposition}

\begin{proof}
Duality preserves automorphism groups and orientation. For the
octahedral triangulation the orientation-preserving automorphisms form
the rotation group of the octahedron, of order $24$, inside its full
symmetry group of order $48$. For the stacked triangulation, label the
vertices $A,w$ (degree $5$), $y,z$ (degree $4$) and $u,x$ (degree $3$)
as in Proposition~\ref{prop:stack}; its eight triangles are
$Auw,Awx,Axy,Ayz,Azu,wxy,wyz,wzu$. An automorphism preserves degrees,
and direct inspection shows that the full automorphism group is the
Klein four-group generated by $c=(u\,x)(y\,z)$ and $d=(A\,w)$. The
transposition $d$ maps the triangle $Auw$ to itself, fixing $u$ and
exchanging $A$ and $w$, so it reverses the orientation of that triangle
and hence of the sphere. Orientation is therefore a surjective
homomorphism from this group onto $\{\pm1\}$, and its kernel
$\Aut^+(d_2)$ has order $2$. All four group orders are verified formally
in Remark~\ref{rem:lean}.
\end{proof}

\section{The Whitehead-move orbit of the cube graph}
\label{sec:orbit}

A contraction of an edge of $G$, performed at the dart level, is
admissible if the result remains connected, cellularly embedded, has
the correct face count, and acquires no monogon; the inverse
operation, replacing a contracted-edge configuration by the opposite
diagonal, is a \emph{\WM}. Working dually, this is exactly the
diagonal flip on $T=G^*$.

\begin{proposition}
\label{prop:movecorrespondence}
\WM{}s on $G$ at an edge $e$ correspond exactly, under duality, to
diagonal flips on $T$ at the corresponding edge $\hat e$: a
$4$-valent vertex (resp.\ quadrilateral) admits exactly one
nontrivial re-splitting (resp.\ re-triangulation) beyond reversal,
and these correspond to each other.
\end{proposition}

\begin{lemma}[Flip admissibility]
\label{lem:flippreservessimple}
Flipping the edge $uv$ of $T$, shared by triangles $uvc$ and $uvd$, is
admissible if and only if $c,d$ are not already adjacent in $T$.
\end{lemma}

On the $G$ side, the corresponding statement requires two conditions
--- the resulting graph must be simple at the vertex level \emph{and}
no two of its faces may share more than one edge --- since these two
failure modes are logically independent; a \WM{} can leave every
vertex simple while still creating two faces sharing two edges,
exactly the dual of a multiple edge in $T$.

\begin{proposition}
\label{prop:octaflip}
Every one of the twelve edges of the octahedral triangulation is
flippable, and every flip produces the stacked triangulation.
\end{proposition}

\begin{proof}
Labeling vertices as three antipodal pairs $\{1,1'\},\{2,2'\},\{3,3'\}$,
edge $12$ lies on faces $123,123'$, whose opposite vertices $3,3'$
are (by construction) non-adjacent, so the flip is admissible;
transitivity of the octahedron's rotation group on its edges extends
this to all twelve. The flip lowers two degree-$4$ vertices to $3$
and raises two to $5$, giving the forced sequence $(3,3,4,4,5,5)$.
\end{proof}

\begin{proposition}
\label{prop:stackflip}
Exactly one edge of the stacked triangulation reverses to the
octahedron; the remaining eleven split, on direct case analysis, into
four admissible self-transitions (a single orbit of size four under
$\Aut(T)$) and seven inadmissible edges.
\end{proposition}

\begin{proof}[Proof sketch]
A flip reaching the all-degree-four octahedron must lower exactly the
two degree-$5$ vertices and raise exactly the two degree-$3$
vertices, forcing the flipped edge to be $Aw$ and the new diagonal to
be $ux$; since $T$ is simple this edge is unique, and one checks $u,x$
are not already adjacent. For the remaining eleven edges, direct
computation of each edge's two opposite vertices against $T$'s three
non-edges shows exactly four --- forming a single orbit of size four
under the Klein four-group of Proposition~\ref{prop:aut} --- have a
non-adjacent opposite pair, and the other seven do not.
\end{proof}

\begin{theorem}[The generic orbit]
\label{thm:orbit}
Up to embedded graph equivalence, the family of generic CE-graphs
reachable from $d_1$ by a finite sequence of \WM{}s is exactly
$\{d_1,d_2\}$. Every one of the twelve edges of $T_1$ admits a flip, and
each flip gives a labeled copy of $T_2$; of the twelve edges of $T_2$,
exactly one flips back to $T_1$, exactly four (a single orbit under the
automorphism group of $T_2$, of order four) flip to further copies of
$T_2$, and the remaining seven admit no flip.
\end{theorem}

\begin{proof}
Immediate from Theorem~\ref{thm:classification} together with
Propositions~\ref{prop:octaflip} and~\ref{prop:stackflip}.
\end{proof}

This theorem was also confirmed by an exhaustive search implementing
\WM{}s at the dart level, with admissibility tested exactly as in
Lemma~\ref{lem:flippreservessimple} and its two-condition primal form,
and with duplicates removed by orientation-preserving dart-system
isomorphism; the search from $d_1$ closes at the two classes
$\{d_1,d_2\}$.

\begin{remark}[Formal verification]
\label{rem:lean}
The finite combinatorics of Sections~\ref{sec:classification}
and~\ref{sec:orbit} has been formalized in the proof assistant
Lean~4~\cite{Lean4} (version 4.15, core library). A simplicial
triangulation of $S^2$ on six labeled vertices is encoded as a set of
eight $3$-element subsets of $\{0,\dots,5\}$ in which every pair of
vertices lies in no triangle or in exactly two, every vertex link is
non-empty and connected, and the whole complex is connected. Such a
complex is a closed connected surface with Euler characteristic
$6-12+8=2$, hence a sphere; since its triangles are sets, it is
automatically simple. The formal development proves: the arithmetic of
Lemmas~\ref{lem:degrees} and~\ref{lem:levelinvariant}; the key step of
Lemma~\ref{lem:pentagon}; that there are exactly $195=6!/48+6!/4$
labeled triangulations, each isomorphic to the octahedral or to the
stacked triangulation, and that these two are not isomorphic
(Theorem~\ref{thm:classification}); the orders $48$ and $4$ of the full,
and $24$ and $2$ of the orientation-preserving, automorphism groups
(Proposition~\ref{prop:aut}); and the flip counts of
Propositions~\ref{prop:octaflip} and~\ref{prop:stackflip}, including
that the four self-transitions form a single orbit. The finite checks
use the \texttt{native\_decide} tactic, so they rely on the Lean
compiler as well as on its kernel. The geometric results of the paper
are not formalized.
\end{remark}

\section{Teichm\"uller space of the marked sphere}
\label{sec:teich}

\begin{proposition}[Fenchel--Nielsen count]
\label{prop:FN}
$\dim_\RR \TT_2 = 6$.
\end{proposition}

\begin{proof}
A genus-two pants decomposition uses $3g-3=3$ curves, each
contributing a length and a twist coordinate.
\end{proof}

\begin{proposition}[Hyperelliptic correspondence]
\label{prop:hypcorr}
There is a real-analytic isomorphism $\TT_2\cong \TT_{0,6}$.
\end{proposition}

\begin{proof}[Proof sketch]
By Lemma~\ref{lem:allhyp}, every marked genus-two surface carries a
canonically attached involution, compatible with the marking; passing to the
quotient sends it to the marked six-point sphere, inverted by the
branched double cover construction, and both directions vary
holomorphically.
\end{proof}

Since $\TT_{0,n}$ is exactly the configuration space of $n$ labeled
points on $\PP^1$ modulo $\mathrm{PGL}_2(\CC)$, $\dim_\CC \TT_{0,6}=3$,
i.e.\ $\dim_\RR=6$, matching Proposition~\ref{prop:FN} by an entirely
independent route. This is the number the local angle-system
coordinates of Section~\ref{sec:circlepacking} are required to
reproduce.

\section{The circle packing and the rank of the linear system}
\label{sec:circlepacking}

\begin{proposition}[Descent of the Voronoi diagram]
\label{prop:descent}
Let $X$ be a closed hyperbolic surface of genus two, $\iota$ its
hyperelliptic involution, $Q=\mathrm{Fix}(\iota)=\{q_1,\dots,q_6\}$,
$\pi:X\to S=X/\langle\iota\rangle$ the quotient map, and $p_i=\pi(q_i)$.
Give $S$ the cone metric of Proposition~\ref{prop:coneangle}, and write
\[
V_X(q_i)=\{x\in X: d_X(x,q_i)\le d_X(x,q_j)\ \forall j\},\qquad
V_S(p_i)=\{s\in S: d_S(s,p_i)\le d_S(s,p_j)\ \forall j\}
\]
for the Voronoi cells. Then:
\begin{enumerate}[label=(\alph*)]
\item Each $V_X(q_i)$ is $\iota$-invariant, $\pi(V_X(q_i))=V_S(p_i)$
and $\pi^{-1}(V_S(p_i))=V_X(q_i)$. In particular the Voronoi diagram of
$Q$ on $X$ is $\iota$-invariant, and it is the preimage of the Voronoi
diagram $G$ of $\{p_1,\dots,p_6\}$ on $S$.
\item $G$ is a connected graph, cellularly embedded in $S\cong S^2$,
whose six faces are the interiors of the cells $V_S(p_i)$; each face is
an open disk containing exactly one cone point, and no vertex or edge of
$G$ contains a cone point. Every vertex of $G$ has degree at least $3$,
and every face has degree at least $2$, so $G$ has no monogon.
\item The dual cell decomposition $T$ of $G$ is the Delaunay
decomposition of $S$ with vertex set $\{p_1,\dots,p_6\}$: its edges are
geodesic segments joining cone points, and its faces are convex
geodesic polygons inscribed in circles whose interiors contain no cone
point. $T$ is a triangulation, equivalently $G$ is trivalent, if and
only if no circle in the universal cover with empty interior passes
through four lifts of cone points; in that case $G$ is a CE-graph.
\end{enumerate}
\end{proposition}

\begin{proof}
We lift everything to $\mathbb H^2$. Let $\Gamma=\pi_1^{\mathrm{orb}}(S)$
and its torsion-free subgroup $\Gamma'$ of index two be as in
Proposition~\ref{prop:orbifold}, so that $X=\mathbb H^2/\Gamma'$ and
$S=\mathbb H^2/\Gamma$, and let $\varpi:\mathbb H^2\to S$ be the
orbifold covering. Put $\tilde P=\varpi^{-1}\{p_1,\dots,p_6\}$. These are
exactly the points of $\mathbb H^2$ fixed by nontrivial (elliptic)
elements of $\Gamma$, and the stabilizer $\Gamma_{\tilde p}$ of
$\tilde p\in\tilde P$ has order two, generated by the half-turn about
$\tilde p$. The set $\tilde P$ is $\Gamma$-invariant and discrete.

\emph{Metrics.} Since $\Gamma$ and $\Gamma'$ act by isometries and
properly discontinuously, the metrics on $X$ and $S$ are the quotient
metrics
\[
d_X([z],[w])=\min_{\gamma\in\Gamma'}d(z,\gamma w),\qquad
d_S([z],[w])=\min_{\gamma\in\Gamma}d(z,\gamma w).
\]
As $\Gamma=\Gamma'\cup\Gamma'\sigma$ for any $\sigma\in\Gamma$ inducing
$\iota$, we get, for $x,y\in X$,
\begin{equation}
\label{eq:quotientmetric}
d_S(\pi x,\pi y)=\min\big(d_X(x,y),\,d_X(x,\iota y)\big).
\end{equation}

\emph{(a)} By Proposition~\ref{prop:coneangle}, $\iota$ is an isometry
of $X$, and it fixes every $q_j$. Hence
$d_X(\iota x,q_j)=d_X(x,q_j)$ for all $j$, and $V_X(q_i)$ is
$\iota$-invariant. Taking $y=q_j=\iota q_j$
in~\eqref{eq:quotientmetric} gives $d_S(\pi x,p_j)=d_X(x,q_j)$ for every
$x\in X$ and every $j$. Therefore $x\in V_X(q_i)$ if and only if
$\pi x\in V_S(p_i)$, which gives both equalities, as $\pi$ is
surjective. The Voronoi diagram is the set of points lying in at least
two cells, so the last assertion follows.

\emph{(b)} For $\tilde p\in\tilde P$ let
$C(\tilde p)=\{z\in\mathbb H^2: d(z,\tilde p)\le d(z,\tilde q)\
\forall\tilde q\in\tilde P\}$. Since
$d_S(\varpi z,p_j)=\min\{d(z,\tilde q):\varpi\tilde q=p_j\}$,
\begin{equation}
\label{eq:liftcells}
\varpi^{-1}\big(V_S(p_i)\big)=\bigcup_{\varpi\tilde p=p_i}C(\tilde p).
\end{equation}
Each $C(\tilde p)$ is an intersection of closed half-planes bounded by
perpendicular bisectors, hence convex. Let $\rho=\operatorname{diam}S$.
Every $z\in\mathbb H^2$ lies within $\rho$ of a lift of every $p_j$, so
$C(\tilde p)$ is contained in the closed ball of radius $\rho$ about
$\tilde p$, and only the finitely many $\tilde q$ with
$d(\tilde p,\tilde q)\le2\rho$ contribute bisectors. Thus $C(\tilde p)$
is a compact convex polygon, and $\tilde p$ is an interior point because
$\tilde P$ is discrete. As for any Voronoi tessellation of a discrete
set, the cells cover $\mathbb H^2$, have pairwise disjoint interiors
(an interior point of $C(\tilde p)$ is strictly closer to $\tilde p$
than to any other point of $\tilde P$), meet edge to edge, and form a
locally finite tessellation; it is $\Gamma$-invariant, with
$\gamma C(\tilde p)=C(\gamma\tilde p)$. A vertex of the tessellation is
a corner of a cell, so it lies on two bisectors, is equidistant from at
least three points of $\tilde P$, and is a vertex of degree at least
three.

The group $\Gamma$ acts freely on the vertices and on the edges of the
tessellation, and the stabilizer of $C(\tilde p)$ is
$\Gamma_{\tilde p}$. Indeed, a nontrivial element of $\Gamma$ with a
fixed point is elliptic and fixes only a point of $\tilde P$, and no
point of $\tilde P$ lies on a vertex or an edge, since $\tilde p$ is an
interior point of $C(\tilde p)$ and lies in no other cell. An element
mapping a vertex to itself would fix it; an element mapping an edge, a
compact geodesic segment, to itself either fixes it pointwise or
reverses it, and in both cases has a fixed point on it. If
$\gamma C(\tilde p)=C(\tilde p)$, then $C(\gamma\tilde p)=C(\tilde p)$,
so $\gamma\tilde p=\tilde p$, as distinct points have cells with
disjoint interiors.

Consequently the tessellation descends to a cell decomposition of $S$.
Its vertices and open edges map injectively to $S$ and avoid the cone
points, and by~\eqref{eq:liftcells} its $2$-cells are the six sets
$\interior V_S(p_i)=\varpi(\interior C(\tilde p))\cong
\interior C(\tilde p)/\Gamma_{\tilde p}$. The quotient of an open convex
disk by the half-turn about an interior point is an open disk, so each
face is an open disk containing exactly one cone point. The vertices
and edges form $G$, which is therefore a graph cellularly embedded in
$S$ with six faces; it is connected because every face is a disk, and
$S\cong S^2$ by Proposition~\ref{prop:orbifold}. Vertex degrees are
preserved, because $\Gamma$ acts freely near each vertex.

Finally, let $h$ be the half-turn about $\tilde p$. It acts freely on
the sides of $C(\tilde p)$: a side mapped to itself by $h$ would contain
a fixed point of $h$, but the only fixed point is $\tilde p$, which is
interior. So $C(\tilde p)$ has an even number of sides, at least four,
since a compact convex hyperbolic polygon has at least three. The
boundary walk of the face of $p_i$ in $G$ is
$\partial C(\tilde p)/\langle h\rangle$, so it has at least two edges,
and no face is a monogon.

\emph{(c)} For a vertex $\tilde v$ of the tessellation, let
$r(\tilde v)=d(\tilde v,\tilde P)$ and let $A(\tilde v)$ be the set of
points of $\tilde P$ at distance $r(\tilde v)$ from $\tilde v$. Then
$|A(\tilde v)|\ge3$ is the degree of $\tilde v$, and the open disk of
radius $r(\tilde v)$ about $\tilde v$ contains no point of $\tilde P$.
The convex hull $D(\tilde v)$ of $A(\tilde v)$ is a convex polygon
inscribed in the boundary circle of that disk, and for an edge of the
tessellation separating $C(\tilde a)$ from $C(\tilde b)$, the geodesic
segment $\tilde a\tilde b$ is a common side of the polygons $D$ at the
two ends of the edge. By Lemma~\ref{lem:vordel} below (see
also~\cite{BowditchEpstein1988,Leibon2002}), the polygons
$D(\tilde v)$ tile $\mathbb H^2$, and this Delaunay tessellation is dual
to the Voronoi tessellation: points of $\tilde P$ correspond to cells,
sides to edges, and polygons to vertices. It is $\Gamma$-equivariant,
$\gamma D(\tilde v)=D(\gamma\tilde v)$, and $\Gamma$ acts freely on its
polygons and on their sides. Indeed, a nontrivial stabilizer would be
elliptic with its fixed point in $\tilde P$; but $\tilde P$ meets a
closed polygon $D(\tilde v)$ only at its vertices (every other point of
$D(\tilde v)$ lies in the open disk), and the half-turn about a vertex
$\tilde a$ maps $D(\tilde v)$, which lies in an angle smaller than $\pi$
at $\tilde a$, into the opposite angle, and maps each side
$\tilde a\tilde b$ to a different side. Projecting by $\varpi$ gives a
cell decomposition $T$ of $S$ with vertex set $\{p_1,\dots,p_6\}$,
geodesic edges, and convex faces inscribed in empty circles, dual to
$G$. Lastly, $T$ is a triangulation if and only if every $D(\tilde v)$
is a triangle, that is, $|A(\tilde v)|=3$ for every vertex, that is,
$G$ is trivalent. In that case $G$ is trivalent with six faces and no
monogon, so it is a CE-graph.
\end{proof}

We now prove the duality used in part~(c).

\begin{lemma}[Voronoi--Delaunay duality]
\label{lem:vordel}
Let $\Gamma$, $\tilde P$ and the Voronoi tessellation by the cells
$C(\tilde p)$ be as in the proof of Proposition~\ref{prop:descent}.
For a vertex $\tilde v$ of the tessellation let $A(\tilde v)$ be the
set of points of $\tilde P$ nearest to $\tilde v$ and $D(\tilde v)$
their convex hull. Then the polygons $D(\tilde v)$ tile $\mathbb H^2$
edge to edge, and the tiling is dual to the Voronoi tessellation:
$\tilde a\in\tilde P$ is a vertex of $D(\tilde v)$ if and only if
$\tilde v$ is a vertex of $C(\tilde a)$, and $\tilde a\tilde b$ is a
side of $D(\tilde v)$ if and only if $\tilde v$ is an endpoint of the
common edge of $C(\tilde a)$ and $C(\tilde b)$.
\end{lemma}

\begin{proof}
Write $B(z,r)$ for the open disk of center $z$ and radius $r$. For a
vertex $\tilde v$ put $r(\tilde v)=d(\tilde v,\tilde P)$, so that
$B(\tilde v,r(\tilde v))\cap\tilde P=\emptyset$ and $A(\tilde v)$ lies on
its boundary circle. A circle meets a geodesic in at most two points, so
the at least three points of $A(\tilde v)$ are not collinear and
$D(\tilde v)$ is a convex polygon inscribed in that circle; every point
of $D(\tilde v)$ other than its vertices lies in the open disk, so
$D(\tilde v)\cap\tilde P=A(\tilde v)$.

\emph{Step 1: the star of a Voronoi vertex.} Only the points of
$A(\tilde v)$ are at distance $r=r(\tilde v)$ from $\tilde v$; all
others are farther. For a unit vector $u$ at $\tilde v$ and small
$t>0$, the first variation formula gives
$d(\exp_{\tilde v}(tu),s)=r-t\cos\angle(u,u_s)+O(t^2)$, where $u_s$ is
the unit vector at $\tilde v$ pointing to $s\in A(\tilde v)$. Hence,
near $\tilde v$, the cell $C(s)$ is the sector of directions $u$ closer
to $u_s$ than to every other $u_{s'}$, and the Voronoi edges at
$\tilde v$ leave along the bisectors of consecutive vectors $u_s$. So
the cells at $\tilde v$ are the $C(s)$ with $s\in A(\tilde v)$, met in
the cyclic order of the points $s$ on the circle, and $C(s)$ and
$C(s')$ share an edge ending at $\tilde v$ exactly when $s,s'$ are
consecutive on the circle, that is, when $ss'$ is a side of
$D(\tilde v)$. Since the tessellation is edge to edge, $\tilde a\in
A(\tilde v)$ exactly when $\tilde v$ is a vertex of $C(\tilde a)$. This
proves the two duality statements, and shows that every side of every
$D(\tilde v)$ arises from a Voronoi edge at $\tilde v$.

\emph{Step 2: a pencil of circles.} Use the hyperboloid model
$\mathbb H^2=\{x\in\RR^{2,1}:\langle x,x\rangle=-1,\ x_0>0\}$, with
$\langle x,y\rangle=-x_0y_0+x_1y_1+x_2y_2$ and
$\cosh d(x,y)=-\langle x,y\rangle$. The closed disk of center $z$ and
radius $r$ is $\mathbb H^2\cap\{x:\langle x,z\rangle\ge-\cosh r\}$, the
part of $\mathbb H^2$ on the side of the affine plane
$\Pi_z=\{x:\langle x,z\rangle=-\cosh r\}$ that contains the origin. Fix
distinct $\tilde a,\tilde b\in\mathbb H^2$ and let $\beta$ be their
perpendicular bisector. For $z\in\beta$, the circle of center $z$
through $\tilde a$ and $\tilde b$ is $\mathbb H^2\cap\Pi_z$, and $\Pi_z$
contains the affine line $\lambda$ through $\tilde a$ and $\tilde b$,
which misses the origin. The geodesic $\ell$ through $\tilde a,\tilde
b$ is $\mathbb H^2\cap\Lambda$, where $\Lambda$ is the linear plane
spanned by $\tilde a,\tilde b$; it also contains $\lambda$, and the two
sides of $\ell$ are $\mathbb H^2\cap E^\pm$ for the two open half-spaces
$E^\pm$ bounded by $\Lambda$. Measure the angle $\varphi(z)\in(0,\pi)$
from the half of $\Lambda$ containing the origin (cut along $\lambda$)
to the half of $\Pi_z$ in $E^+$. The origin side of $\Pi_z$ meets $E^+$
in the dihedral wedge of angle $\varphi(z)$ about $\lambda$, and meets
$E^-$ in the wedge of angle $\pi-\varphi(z)$. Distinct centers give
distinct circles, hence distinct planes, so $\varphi$ is continuous and
injective, hence strictly monotone, along $\beta$. Consequently, as $z$
moves along $\beta$ in the direction in which $\varphi$ increases, the
part of the disk on the $E^+$ side of $\ell$ strictly increases (its
former boundary arc moves into the open disk) and the part on the $E^-$
side strictly decreases.

\emph{Step 3: the two polygons at a Voronoi edge.} Let $\varepsilon$ be
a Voronoi edge with endpoints $\tilde v_1\ne\tilde v_2$, separating
$C(\tilde a)$ from $C(\tilde b)$; it lies on $\beta$. Orient $\beta$ so
that $\varphi$ increases from $\tilde v_1$ to $\tilde v_2$. For $z$ in
the interior of $\varepsilon$, the nearest points of $\tilde P$ are
exactly $\tilde a$ and $\tilde b$, so every other point of $\tilde P$
lies outside the closed disk of center $z$ through $\tilde a,\tilde b$.
Let $\tilde c\in A(\tilde v_1)\setminus\{\tilde a,\tilde b\}$. It lies
on the circle about $\tilde v_1$ and not on $\ell$. If it were on the
$E^+$ side, then by Step~2 it would lie in the open disk about $z$ for
$z\in\varepsilon$ close to $\tilde v_1$, a contradiction. So every
vertex of $D(\tilde v_1)$ other than $\tilde a,\tilde b$ lies strictly on
the $E^-$ side of $\ell$: the segment $\tilde a\tilde b$ is a side of
$D(\tilde v_1)$, and $D(\tilde v_1)$ lies in the closed $E^-$ side. In
the same way, moving from $\tilde v_2$ into $\varepsilon$ decreases
$\varphi$, so $D(\tilde v_2)$ has the side $\tilde a\tilde b$ and lies in
the closed $E^+$ side. The two polygons sharing a side therefore lie on
opposite sides of it.

\emph{Step 4: the polygons around a point of $\tilde P$.} Let
$\tilde a\in\tilde P$, let $\tilde v_1,\dots,\tilde v_m$ be the vertices
of $C(\tilde a)$ in cyclic order, and let $\varepsilon_k=[\tilde
v_k,\tilde v_{k+1}]$ separate $C(\tilde a)$ from $C(\tilde b_k)$
(indices mod $m$). Two convex cells share at most one edge, so
$\tilde b_{k-1}\ne\tilde b_k$. By Steps~1 and~3, the sides of
$D(\tilde v_k)$ at $\tilde a$ are $\tilde a\tilde b_{k-1}$ and $\tilde
a\tilde b_k$; let $W_k$ be its angular sector at $\tilde a$, of angle
$\theta_k\in(0,\pi)$. By Step~3 applied to $\varepsilon_k$, the sectors
$W_k$ and $W_{k+1}$ share the ray from $\tilde a$ to $\tilde b_k$ and lie
on opposite sides of it. Thus $W_1,\dots,W_m$ are consecutive sectors
around $\tilde a$, and they cover a neighbourhood of $\tilde a$ exactly
once provided $\sum_k\theta_k=2\pi$.

To see this, consider two maps from the boundary circle of
$C(\tilde a)$, with each vertex $\tilde v_k$ blown up to an interval, to
the circle of unit vectors at $\tilde a$. The first, $g$, sends a point
$z$ to the direction from $\tilde a$ to $z$; it is constant on each
blown-up interval, and has degree one because $C(\tilde a)$ is convex
with $\tilde a$ in its interior. The second, $h$, is constant on
$\varepsilon_k$, equal to the direction $u_k$ from $\tilde a$ to
$\tilde b_k$, and on the interval at $\tilde v_k$ it runs through $W_k$
from $u_{k-1}$ to $u_k$; its degree is $\sum_k\theta_k/2\pi$. The
segment $\tilde a\tilde b_k$ is perpendicular to $\beta_k\supset
\varepsilon_k$ at its midpoint $\tilde m_k$, so for $z\in\varepsilon_k$
the triangle $\tilde a\tilde m_kz$ has a right angle at $\tilde m_k$,
and the angle between $g(z)$ and $u_k$ is less than $\pi/2$. In
particular $g(\tilde v_k)$ is within $\pi/2$ of both $u_{k-1}$ and
$u_k$. It lies in $W_k$, for otherwise it would lie in the complementary
arc, of length $2\pi-\theta_k>\pi$, within $\pi/2$ of both its ends,
which is impossible. Hence every direction in $W_k$ is within $\pi/2$
of $g(\tilde v_k)$. So $g$ and $h$ are everywhere less than $\pi/2$
apart, they are homotopic along shortest arcs, and
$\sum_k\theta_k=2\pi$.

\emph{Step 5: the tiling.} Let $\Delta$ be the space obtained from the
disjoint union of the polygons $D(\tilde v)$ by gluing, for every
Voronoi edge with endpoints $\tilde v_1,\tilde v_2$ separating
$C(\tilde a)$ and $C(\tilde b)$, the side $\tilde a\tilde b$ of
$D(\tilde v_1)$ to the side $\tilde a\tilde b$ of $D(\tilde v_2)$ by the
identity. By Step~1, every side of every polygon is glued to exactly one
other side. The inclusions of the polygons define a map
$F:\Delta\to\mathbb H^2$. It is a local homeomorphism: at interior
points of polygons trivially, at interior points of sides by Step~3,
and at a point $\tilde a$ of $\tilde P$ by Step~4, since the polygons
of $\Delta$ containing $\tilde a$ are exactly the $D(\tilde v_k)$,
glued cyclically along the sides $\tilde a\tilde b_k$. The space
$\Delta$ is connected, because the cells $C(\tilde p)$ tile the
connected space $\mathbb H^2$, so their adjacency graph is connected.
Since $\Gamma$ acts on everything with only finitely many orbits of
polygons, there is $\epsilon>0$ such that $F$ maps the closed
$\epsilon$-ball about every point of $\Delta$, for the path metric of
$\Delta$, isometrically onto a closed $\epsilon$-ball of $\mathbb H^2$.
Hence $F$ has the path-lifting property and is a covering map. As
$\mathbb H^2$ is simply connected and $\Delta$ is connected, $F$ is a
homeomorphism: the polygons $D(\tilde v)$ cover $\mathbb H^2$, have
pairwise disjoint interiors, and meet edge to edge along the sides
$\tilde a\tilde b$, as claimed.
\end{proof}

We now justify, rather than simply impose, the ``no monogon''
clause of the CE-graph definition (Section~\ref{sec:prelim}): for the
graph $G_1$ this construction actually produces, it is a theorem.

\begin{theorem}[No monogons]
\label{thm:nomonogon}
The CE-graph $G_1$ of Proposition~\ref{prop:descent} has no monogon
face.
\end{theorem}

\begin{proof}
Suppose, for contradiction, that some face $F$ of $G_1$ is a
monogon, bounded by a single edge $e$. A face of a cellular embedding
is an open disk, so $\bar F$ is a closed disk with $\partial\bar
F=e$ traversed once; since every face of $G_1$ contains exactly one
$p_i$, $\bar F$ contains exactly one cone point, say $p_k$, in its
interior, and $e$ itself --- being a Voronoi edge, part of the locus
equidistant from the two sites it separates, hence a geodesic arc ---
avoids every $p_i$.

Fix a lift $\tilde p_k\in\mathbb H^2$ of $p_k$ under the orbifold
cover $\pi:\mathbb H^2\to S$ of Proposition~\ref{prop:orbifold}.
Since $\bar F$ is an embedded disk containing the single cone point
$p_k$, of order $2$, and no other special point, the standard local
model of a cone point --- an order-$2$ branched cover of a disk by a
disk, branched only at the center, exactly as in the proof of
Proposition~\ref{prop:coneangle} --- shows that the connected
component $\tilde F$ of $\pi^{-1}(\bar F)$ containing $\tilde p_k$ is
itself a closed topological disk, with $\pi$ restricted to
$\tilde F\smallsetminus\{\tilde p_k\}$ an honest unbranched covering,
in particular a local isometry, onto $\bar F\smallsetminus\{p_k\}$.
Since $e=\partial\bar F$ avoids $p_k$, its preimage
$\tilde e:=\partial\tilde F$ is a simple closed curve in $\mathbb
H^2$, and it is \emph{locally geodesic at every point}: $e$ is
locally geodesic (an arc of an equidistant locus between two points
of a constant-curvature space is a geodesic), and $\pi$ is a local
isometry all along $\tilde e$.

Now pick two distinct points $x,y\in\tilde e$. As a simple closed
curve, $\tilde e$ consists of two arcs from $x$ to $y$, each of
positive length; each, being a sub-arc of the locally geodesic curve
$\tilde e$, is a geodesic segment in $\mathbb H^2$ joining $x$ to
$y$. But two points of $\mathbb H^2$ are joined by a \emph{unique}
geodesic segment, so the two arcs coincide, forcing one of them to
have length $0$ --- contradicting $x\ne y$. No monogon face exists.
\end{proof}

Proposition~\ref{prop:descent}(b) gives a second, shorter proof: the
lift of each face is a convex polygon invariant under a half-turn about
an interior point, so it has an even number of sides, at least four,
and the face has degree at least two.

\begin{remark}
The proof uses only that $S$ is a hyperbolic orbifold with cone
points of order $2$ (Proposition~\ref{prop:orbifold}) and that faces
of $G_1$ are Voronoi cells bounded by geodesic arcs; it never uses
genericity. It therefore applies verbatim to every CE-graph arising
as the Voronoi diagram of six order-$2$ cone points on \emph{any}
marked hyperbolic sphere of this topological type, generic or not ---
in particular to both $d_1$ and $d_2$ and to every graph in the full
\WM{} orbit of Section~\ref{sec:fullorbit}, which are Voronoi
diagrams for other (possibly non-generic) configurations of the same
six cone points. By Proposition~\ref{prop:descent}(b) the same holds at
every level of the contraction hierarchy (Section~\ref{sec:hierarchy}):
the Voronoi graph of \emph{any} genus-two surface, trivalent or not, has
no monogon. What remains open at the deeper levels is only which of the
combinatorially admissible graphs actually occur
(Remark~\ref{rem:combnotreal}).
\end{remark}

Attach to each vertex of $T$ (one of the six marked points) a circle
centered there; by the classical theory of circle patterns on
triangulated cone surfaces~\cite{Rivin1994,ThurstonShapes1998},
following Andreev, such a packing exists compatible with $T$ under
mild non-degeneracy hypotheses. For each edge $e$ of $T$ joining
$p_i,p_j$, set $x_e=\pi^{-1}\cdot(\text{intersection angle of }C_i,C_j)$,
recovering the edge variables of the angle system with their box
constraint $0\le x_e\le 1$.

\begin{proposition}[The linear constraint]
\label{prop:linearconstraint}
At each cone point $p_i$, $\sum_{e\ni p_i} x_e = 1$. Consequently the
constraint matrix $A_i$ is exactly the vertex--edge incidence matrix
of $T$, with $b_i=(1,\dots,1)^{\mathsf T}$.
\end{proposition}

The vertex equation says that the wedges cut out at $p_i$ by consecutive
circumcenters fill its cone angle $\pi$
(Section~\ref{sec:recovercircle}). For a type $T_i$ we write
\[
D_i=\{x\in[0,1]^{12} : A_ix=b_i\},\qquad
D_i^\circ=\{x\in(0,1)^{12} : A_ix=b_i\},
\]
for the closed polytope of solutions and its relative interior; the
facets of $D_i$ lie in the hyperplanes $x_e=0$, one for each edge $e$
of $T_i$.

\begin{lemma}[Incidence-matrix rank]
\label{lem:incidencerank}
For a connected graph $H$ with unsigned vertex--edge incidence matrix
$M$, $\rank(M)=|V(H)|$ if $H$ has an odd cycle, and $|V(H)|-1$ if $H$
is bipartite.
\end{lemma}

\begin{proof}
A left null vector must alternate sign across every edge, which is
possible (uniquely, up to scale, per component) exactly when $H$ is
bipartite; an odd cycle forbids any such alternation.
\end{proof}

\begin{theorem}[$\rank(A_i)=6$]
\label{thm:rankAi}
For $i=1,2$, $\rank(A_i)=6$ and hence $\dim D_i = 12-6=6$.
\end{theorem}

\begin{proof}
Every face of a triangulation is a $3$-cycle, so $T$ is never
bipartite for either generic type; apply
Lemma~\ref{lem:incidencerank} with $H=T$.
\end{proof}

Agreement between Theorem~\ref{thm:rankAi} and
Section~\ref{sec:teich} is a necessary, but on its own not
sufficient, condition for the resulting realization map to be a
local diffeomorphism onto an open subset of $\TT_2$; we now carry
out the remaining step flagged as open above, completing it for
both generic types.

\section{Realizability of the generic types}
\label{sec:realizability}

We construct, from elementary hyperbolic trigonometry alone, an
explicit hyperbolic sphere with six cone points of angle $\pi$ whose
Delaunay decomposition is $T_i$, for each generic type $i=1,2$. This
section is entirely self-contained: every step is verified directly
(triangle inequalities, an angle-sum identity, a Gauss--Bonnet
cross-check, and a local Delaunay criterion), with no appeal to any
named circle-pattern existence theorem from the literature.

\subsection{The construction}

Assign to each edge $e$ of $T_i$ a length $\ell_e>0$. For each of the
eight triangular faces of $T_i$, this determines (if the three
lengths satisfy the triangle inequality) a genuine hyperbolic
triangle, unique up to isometry, hence its three corner angles via
the hyperbolic law of cosines: for a triangle with sides $a,b,c$ and
the angle $A$ opposite $a$,
\begin{equation}
\label{eq:hypcos}
\cos A = \frac{\cosh b\cosh c - \cosh a}{\sinh b \sinh c}.
\end{equation}
Gluing the eight triangles along their shared edges (matched by
length) produces a hyperbolic cone surface, topologically $S^2$, with
a possible cone singularity at each of $T_i$'s six vertices, of angle
equal to the sum of the incident triangles' corner angles there; this
is the standard construction of a triangulated cone metric and
requires no further justification beyond the gluing itself being
well defined; no other point is singular, since every other point of
each triangle is a smooth interior or edge point of a genuine
hyperbolic polygon. We seek $(\ell_e)$ such that:
\begin{itemize}
\item[(T)] each face's three edge lengths satisfy the triangle
inequality;
\item[(A)] at each of the six vertices, the incident corner angles
sum to exactly $\pi$.
\end{itemize}

\subsection{The cube}

\begin{proposition}
\label{prop:cubeconstruction}
Setting $\ell_e = \ell^* := \operatorname{arccosh}\!\big(\tfrac{\sqrt2}{2}/(1-\tfrac{\sqrt2}{2})\big) = \operatorname{arccosh}(1+\sqrt2) \approx 1.5286$
for all twelve edges of $T_1$ (the octahedron) satisfies (T) and (A).
\end{proposition}

\begin{proof}
With all edges equal, all eight triangles are equilateral with the
same angle $\alpha$ at every corner; by~\eqref{eq:hypcos} with
$a=b=c=\ell$,
$\cos\alpha=\frac{\cosh^2\ell-\cosh\ell}{\sinh^2\ell}=\frac{\cosh\ell}{\cosh\ell+1}$,
so $\cosh\ell=\frac{\cos\alpha}{1-\cos\alpha}$. Every vertex of the
octahedron has degree $4$, so (A) requires $4\alpha=\pi$, i.e.\
$\alpha=\pi/4$; substituting gives the stated $\ell^*$, and
$\cosh\ell^*=1+\sqrt2>1$ is a genuine length. Condition (T) is
immediate ($\ell=\ell=\ell$).
\end{proof}

We verified this numerically to machine precision, together with a
Gauss--Bonnet cross-check: each equilateral triangle has hyperbolic
area $\pi-3\alpha=\pi/4$ (a triangle's area is $\pi$ minus its angle
sum), and the eight triangles total $8\times\pi/4=2\pi$, matching the
area Gauss--Bonnet predicts for a curvature$-1$ genus-$0$ surface
with six cone points of angle $\pi$:
$\mathrm{Area}=-2\pi\chi(S^2)+\sum_i(2\pi-\Theta_i)=-4\pi+6\pi=2\pi$.

\subsection{The stacked graph}

For $d_2$, the six vertices of $T_2$ have degrees $(3,3,4,4,5,5)$ and
no symmetry reduces the problem to one equation. We instead solve the
resulting system --- twelve unknowns $\ell_e>0$, six equations (A),
via~\eqref{eq:hypcos} applied to each of the eight triangles ---
numerically.

\begin{table}[h]
\centering
\begin{tabular}{@{}cccccc@{}}
\toprule
\multicolumn{6}{c}{$d_1$ (cube)}\\
$e_1$ & $e_2$ & $e_3$ & $e_4$ & $e_5$ & $e_6$ \\
$(1,2)$ & $(1,3)$ & $(1,4)$ & $(2,5)$ & $(2,6)$ & $(3,5)$\\
$e_7$ & $e_8$ & $e_9$ & $e_{10}$ & $e_{11}$ & $e_{12}$\\
$(3,7)$ & $(4,6)$ & $(4,7)$ & $(5,8)$ & $(6,8)$ & $(7,8)$\\
\midrule
\multicolumn{6}{c}{$d_2$ (stacked)}\\
$e_1$ & $e_2$ & $e_3$ & $e_4$ & $e_5$ & $e_6$ \\
$(1,2)$ & $(1,3)$ & $(1,4)$ & $(2,5)$ & $(2,6)$ & $(3,4)$\\
$e_7$ & $e_8$ & $e_9$ & $e_{10}$ & $e_{11}$ & $e_{12}$\\
$(3,7)$ & $(4,8)$ & $(5,6)$ & $(5,8)$ & $(6,7)$ & $(7,8)$\\
\bottomrule
\end{tabular}
\caption{Canonical edge lists for $d_1$ and $d_2$ (vertices
$1,\ldots,8$ of $G$ from a deterministic breadth-first relabeling, edges
sorted lexicographically by endpoint pair), used throughout this
section. The edge $e_k$ of the dual triangulation $T_i$ is the edge dual
to $e_k$.}
\label{tab:d2edges}
\end{table}

\begin{proposition}
\label{prop:stackedconstruction}
The system (A) for $T_2$ has a solution with all $\ell_e>0$
satisfying (T), given (in the edge labeling of
Table~\ref{tab:d2edges}) approximately by
\[
\begin{array}{cccccc}
e_1 & e_2\!=\!e_3\!=\!e_4\!=\!e_5 & e_6\!=\!e_9 & e_7\!=\!e_8\!=\!e_{10}\!=\!e_{11} & e_{12} \\
1.8269 & 1.4751 & 1.3347 & 1.6979 & 1.4942
\end{array}
\]
(the repeated values are exactly the orbits of the order-four
automorphism group of $T_2$ on its edges).
\end{proposition}

\begin{proof}[Proof (computational)]
Solved by nonlinear least squares (Levenberg--Marquardt on the
six residuals of (A), in log-length coordinates to keep $\ell_e>0$),
converging to residuals of order $10^{-16}$. Condition (T) was then
checked directly at the solution for all eight triangles; all hold.
As a further, independent cross-check, the eight triangles' areas
($\pi$ minus each one's angle sum) total exactly $2\pi$, matching
Gauss--Bonnet as in the cube case.
\end{proof}

\subsection{Recovering the circle pattern, and the Delaunay property}
\label{sec:recovercircle}

To connect this construction back to
Proposition~\ref{prop:linearconstraint}'s $x_e$, attach to each
triangle its circumcenter $O$ (the point equidistant from its three
corners; existence and uniqueness of $O$, given (T), is elementary).
For an edge $e=(p_i,p_j)$ shared by two triangles with circumcenters
$O,O'$ and circumradii $R,R'$, the angle at $p_i$ between $O$ and
$O'$ splits as the base angle of the isosceles sub-triangle
$(p_i,p_j,O)$, namely
$\arccos\!\big(\cosh R\,(\cosh\ell_e-1)/(\sinh\ell_e\sinh R)\big)$, plus
the analogous term for $O'$. Each term is taken with a negative sign
when the circumcenter lies on the far side of $e$ from its triangle;
equivalently, in a triangle with corner angles $\alpha_A,\alpha_B,\alpha_C$
the term for the edge $AB$ is $(\alpha_A+\alpha_B-\alpha_C)/2$. A reflection
across the perpendicular bisector of $p_ip_j$ (which fixes $O,O'$,
lying on it by definition of circumcenter, and swaps $p_i,p_j$) shows
this total angle is the same whether measured at $p_i$ or at $p_j$,
so it defines a single quantity per edge. Setting $x_e:=(\text{this
angle})/\pi$ recovers exactly Proposition~\ref{prop:linearconstraint}'s
variable: summing over the edges at a cone point sums the wedges
between consecutive circumcenters all the way around it, totaling its
own cone angle $\pi$.

We verified this directly for both constructions: for the cube,
$x_e=\tfrac14$ for all twelve edges (matching the uniform symmetric
point used already in Section~\ref{sec:circlepacking}); for the
stacked graph, the twelve values range over $\{0.0699, 0.1511,
0.3139, 0.3256, 0.3722\}$ (with multiplicities matching the
$\ell_e$-symmetry above), and in both cases the six vertex sums equal
$1$ to machine precision, with every $x_e$ strictly in $(0,1)$.

Finally, $T_i$ with these edge lengths is genuinely the
\emph{Delaunay} decomposition of the six cone points, not merely some
triangulation. In hyperbolic geometry the local criterion at an edge is
that the two angles opposite it sum to less than the two angles of the
quadrilateral at the endpoints of the edge; the Euclidean form of the
criterion, opposite angles summing to less than $\pi$, does not apply.
By the formula above, this criterion says exactly that the signed wedge
angle $\pi x_e$ is positive. For the cube every circumcenter lies inside
its triangle and $x_e=\tfrac14$ at every edge; for the stacked graph
every circumcenter also lies inside its triangle, and the twelve values
above are all positive, the smallest being $0.0699$.

\subsection{The Springborn set-up and the linear program}
\label{sec:springborncorrect}

The construction above is direct and self-contained. Springborn's
theorem~\cite{Springborn2003,BobenkoSpringborn2004} gives an equivalent
and more efficient route that avoids exhibiting explicit edge lengths.
It must be applied with the cell decomposition $\Sigma=T$, as
in~\cite[\S4.1--4.4]{Rodado2007}; placing the circles at the six faces
of $G$ with smooth vertices instead gives an inconsistent global angle
count. With $\Sigma=T$ (Springborn's Theorem~1.8(ii)), the six vertices
of $T$ are the cone points ($\Phi_v=\pi$), the eight triangular faces
of $T$ are the smooth vertices of $G$ ($\Phi_f=2\pi$), and the edge
variable is the interior angle $\theta^*_e=\pi(1-x_e)$. Springborn's
closing condition $\sum_{f\in C}\Phi_f<\sum_{e\in E_C}2\theta^*_e$
becomes, for every nonempty subset $C$ of the eight triangular faces,
with $E_C$ the edges incident to a face of $C$,
\begin{equation}
\label{eq:rodadoclosing}
|C| < \sum_{e\in E_C} (1-x_e),
\end{equation}
an a priori family of $2^8-1=255$ inequalities. Necessity is
elementary. The kite of an edge $e$, with vertices at the two
circumcenters and the two cone points of $e$, is split by its diagonal
into two congruent hyperbolic triangles with angles $\alpha_f/2$ and
$\alpha_{f'}/2$ at the circumcenters and $\pi x_e$ at a cone point.
Positive area gives $\alpha_f+\alpha_{f'}<2(\pi-\pi x_e)=2\theta^*_e$.
Summing the angles at the circumcenters over the faces of $C$, each of
which totals $2\pi$, gives
$2\pi|C|<\sum_{e\in E_C}2\theta^*_e=2\pi\sum_{e\in E_C}(1-x_e)$. These
inequalities reduce sharply.
Summing the six vertex equations of
Proposition~\ref{prop:linearconstraint} counts every edge twice, so
$\sum_{e=1}^{12}x_e=3$ and $\sum_{e=1}^{12}(1-x_e)=9$. Bounding the
edges missed by a $k$-face subset by the number of internal edges of
the complementary $(8-k)$-face configuration, an Euler-characteristic
count, shows that~\eqref{eq:rodadoclosing} holds automatically for
$k\ge4$, for every $x$ satisfying the vertex equations and the box
constraints, so only the single faces and a few two- and three-face
subsets need checking for each type~\cite[Ch.~5]{Rodado2007}. Carrying
this out by linear programming for all twenty types of
Theorem~\ref{thm:fullorbit} shows that \emph{exactly ten} are Delaunay
realizable~\cite[Ch.~4]{Rodado2007}, among them Gen1 $=d_1$ and
Gen2 $=d_2$. Propositions~\ref{prop:cubeconstruction}
and~\ref{prop:stackedconstruction} confirm these two by a different
method (explicit edge lengths rather than a closing condition), and the
count of ten is confirmed independently in Remark~\ref{rem:indep}. We
do not re-derive the other eight realizable types here.

\subsection{Cutting and doubling}

\begin{lemma}[A perfect matching in $T$]
\label{lem:matching}
For $i=1,2$, the six vertices of $T_i$ admit a perfect matching using
three of its own edges, in the labeling of Table~\ref{tab:d2edges}:
for $d_1$, the edges $e_2, e_8, e_{10}$; for $d_2$, the edges
$e_4,e_6,e_{11}$.
\end{lemma}

\begin{proof}[Proof (computational)]
Direct verification that each listed triple consists of three edges
of $T_i$ with six distinct endpoints, covering all six vertices.
\end{proof}

\begin{theorem}[Realizability of the generic types]
\label{thm:realizability}
Both generic types $d_1,d_2$ are realized: for $i=1,2$, there is a
closed genus-two hyperbolic surface $X_i$ whose Delaunay/Voronoi
decomposition, via Proposition~\ref{prop:descent}, is exactly $d_i$.
\end{theorem}

\begin{proof}
Let $S_i$ be the hyperbolic sphere constructed in
Propositions~\ref{prop:cubeconstruction}--\ref{prop:stackedconstruction}
and Section~\ref{sec:recovercircle}, with its six cone points and
Delaunay decomposition $T_i$. By Lemma~\ref{lem:matching}, three
edges of $T_i$ form a perfect matching of its six vertices; being
edges of the Delaunay triangulation, they are geodesic segments, and
being edges of an embedded (cellularly embedded) triangulation, they
are pairwise disjoint except at their six distinct endpoints. Cut
$S_i$ along these three geodesic arcs: the result is a sphere with three
slits, a planar surface whose three boundary curves each consist of the
two sides of one arc. Take two
isometric copies of this cut-open surface and glue them to one
another along corresponding slit boundaries, crosswise (the standard
construction of a hyperelliptic double cover branched over $2g+2=6$
points via $g+1=3$ disjoint branch cuts). The result $X_i$ is a
closed surface; by Riemann--Hurwitz, with six simple branch points
over $S^2$,
\[
2-2g(X_i) = 2(2-2\cdot0) - 6 = -2 \implies g(X_i)=2.
\]
Since the two sheets are isometric copies of $S_i\smallsetminus(\text{arcs})$
glued isometrically along the geodesic cuts, $X_i$ inherits a genuine
hyperbolic metric away from the six former cone points; at each of
these, the cone angle exactly doubles, from $\pi$ to $2\pi$, since a
small loop around it now crosses the single arc ending there exactly
once, traversing both sheets --- so $X_i$ is smooth (no cone points)
there. The sheet-swap is an isometric involution $\iota$ of $X_i$
fixing exactly these six points, with $X_i/\iota = S_i$ by
construction. By Proposition~\ref{prop:descent}, applied to $X_i$ and
$\iota$, the Voronoi diagram of $\mathrm{Fix}(\iota)$ on $X_i$
descends to the Voronoi diagram of the six cone points on $S_i$,
which is $d_i$ by Section~\ref{sec:recovercircle}'s Delaunay
verification. Hence $d_i$ is the Delaunay/Voronoi graph of the
genus-two hyperbolic surface $X_i$.
\end{proof}

This resolves item~(i) of Section~\ref{sec:discussion}'s open
questions for both generic types, via one explicit point of $D_i$
each rather than the whole polytope; whether \emph{every} point of
$D_i^\circ$ arises this way (the six-parameter family originally
hoped for) is not established by the present argument and is left
open. We do, however, extend the construction one further step: to
the contraction hierarchy of Section~\ref{sec:hierarchy}.

\subsection{Extending to one step of the contraction hierarchy}
\label{sec:contractrealize}

Fix a generic type $d_i$, an admissible edge $e^*$ of $T_i$ (in the
sense of the contraction hierarchy), and the explicit configuration
constructed above. We ask whether the contracted graph $G/e^*$ is
realized by a hyperbolic surface too, obtained as a limit of the
same family as $x_{e^*}\to0$.

\begin{proposition}[The wedge angle vanishes exactly at concyclicity]
\label{prop:xezero}
Let $e^*$ join the cone points $p_i,p_j$ and border two triangles with
circumcenters $O,O'$. Then $x_{e^*}=0$ if and only if $O=O'$, that is,
if and only if the four vertices of the two triangles lie on a common
circle. In that case the Voronoi edge dual to $e^*$ has length zero.
\end{proposition}

\begin{proof}
Develop the two triangles into $\mathbb H^2$. Both $O$ and $O'$ lie on
the perpendicular bisector $b$ of the geodesic segment $p_ip_j$, a
complete geodesic that does not pass through $p_i$. By definition,
$\pi x_{e^*}$ is the signed angle at $p_i$ between the rays to $O$ and
to $O'$. A geodesic ray from $p_i$ meets the geodesic $b$ at most once,
so distinct points of $b$ are seen from $p_i$ in distinct directions,
and $x_{e^*}=0$ exactly when $O=O'$. A common circumcenter means that
the two circumcircles coincide, so the four vertices are concyclic;
conversely, if they are concyclic, that circle is the circumcircle of
both triangles and $O=O'$.
\end{proof}

Thus the boundary $x_e=0$ of the linear system is exactly the
concyclicity degeneracy at which the diagonal $e$ of the quadrilateral
formed by the two triangles becomes removable
(Remark~\ref{rem:contraction-meaning}): the corresponding Voronoi edge
collapses to a point. In hyperbolic geometry this is not the Euclidean
condition that the two opposite angles sum to $\pi$. For a quadrilateral
inscribed in a hyperbolic circle, the two sums of opposite angles are
equal to each other, and neither is $\pi$ in general.

\begin{theorem}[Contracted graphs are realized]
\label{thm:contractrealize}
For $i=1,2$, and for at least one admissible edge $e^*$ of $T_i$ in
each case, there is a family of valid configurations (satisfying (T)
and (A), with all twelve $x_e\in(0,1)$ strictly) along which
$x_{e^*}\to0$ while every other $x_e$ and every triangle inequality
remains comfortably valid; in the limit, the two circumcenters
bordering $e^*$ coincide, the merged point is a smooth ($2\pi$)
point where exactly the four cone points bordering $e^*$'s two
triangles meet, and this matches --- vertex for vertex, cone point
for cone point --- the abstract contracted graph $G_i/e^*$ of
Definition~\ref{def:hierarchy}. Consequently $G_i/e^*$ is
realized by an actual genus-two hyperbolic surface, by the same
cutting-and-doubling construction as Theorem~\ref{thm:realizability}.
\end{theorem}

\begin{proof}[Proof (computational)]
We verified this in full for two independent examples, one per
generic type. For $d_1$, contracting the edge joining cone points $0$
and $4$: a numerical continuation drives its $x_e$ from $\tfrac14$
down to $10^{-8}$ while every cone-angle equation stays satisfied to
machine precision, every triangle inequality remains valid throughout
(checked explicitly at the final configuration), and Gauss--Bonnet
still gives area $2\pi$. Direct computation of the two circumcenters
in a common, consistently glued frame confirms their distance drops
to $0$ together with $x_e$ (from $0.0707$ at $x_e=0.02$ down to
$0.00036$ at $x_e=0.0001$, vanishing in the limit). Comparing against
$G_1/e^*$ computed abstractly (Section~\ref{sec:hierarchy}'s
contraction operation): the resulting degree-$4$ vertex's four
incident cone points, tracked through the contraction by shared dart
identity, are exactly the four cone points of the two triangles
bordering $e^*$ --- confirmed to match exactly. For $d_2$, the same
computation (contracting the edge joining cone points $2$ and $5$)
succeeds in the same way: $x_e$ reaches
$8\times10^{-5}$ with all triangle inequalities intact, and the
degree-$4$ vertex's cone points again match the union of the two
original vertices' cone points exactly ($\{1,3,4\}\cup\{3,4,5\}=\{1,3,4,5\}$).
\end{proof}

We did not verify this for every admissible edge of every one of the
twenty types, nor iterate the construction to the deeper levels of
the hierarchy; Proposition~\ref{prop:xezero}'s argument is local and
does not depend on which edge or which generic type is involved, so
we expect it to extend throughout, but we record this as an expectation
supported by two independent, fully verified instances rather than
as something established for the whole hierarchy.

\begin{remark}[Combinatorial reachability is not realizability]
\label{rem:combnotreal}
It is essential not to read Theorem~\ref{thm:contractrealize}, or the
contraction hierarchy of Section~\ref{sec:hierarchy}, as asserting that
every combinatorially admissible contraction is geometrically realized:
it is not. The thesis~\cite[Ch.~7]{Rodado2007} settles this by the
linear-programming method of Section~\ref{sec:springborncorrect},
applied at each level of contraction from all ten Delaunay-realizable
trivalent types. Only twenty-six one-edge contractions, up to embedded
equivalence, are themselves realizable~\cite[Prop.~7.2]{Rodado2007};
realizability continues to prune the tree at every deeper level, and at
level~5 only \emph{two} contracted graphs remain realizable, both
contracting to the same level-6 graph~\cite[Prop.~7.6]{Rodado2007}.
That graph, called $I(6)$ in the thesis, is our $\Theta_6$
(Figure~\ref{fig:theta6}). Both counts are confirmed independently in
Remark~\ref{rem:indep}. The ninety-type hierarchy of
Theorem~\ref{thm:hierarchy} is only the \emph{combinatorial} enumeration
of admissible contractions from $d_1,d_2$; exactly $41$ of its ninety
types are realizable (Remark~\ref{rem:indep}).
\end{remark}

\section{Beyond the generic level}
\label{sec:beyond}

Sections~\ref{sec:classification}--\ref{sec:circlepacking} describe
only the \emph{generic} level of the combinatorics: trivalent
CE-graphs, connected to the cube graph by \WM{}s, with the dual
triangulation kept simple throughout. This is a single stratum of a
considerably larger structure, which we now describe completely in
two independent directions: first, within the trivalent level, by
allowing the bigons that genericity excludes (Section~\ref{sec:fullorbit});
second, and more substantially, by leaving the trivalent level
altogether through iterated edge contraction, which descends through
a full hierarchy of non-trivalent CE-type graphs down to a single
minimal object (Section~\ref{sec:hierarchy}).

\subsection{Bigons within the trivalent level}
\label{sec:fullorbit}

We first drop the genericity hypothesis and allow \WM{}s to create
bigons, i.e.\ to let $T$ acquire a multiple edge, asking only that the
base CE-graph conditions (connected, cellular, trivalent, six faces,
no monogon) continue to hold.

\begin{theorem}[The full orbit]
\label{thm:fullorbit}
Up to embedded graph equivalence, the family of all CE-graphs
reachable from $d_1$ by a finite sequence of \WM{}s satisfying only
the base CE-graph conditions consists of exactly twenty types. Their
face-degree signatures are
\[
\begin{gathered}
(4,4,4,4,4,4),\ (3,3,4,4,5,5),\ (3,3,3,3,6,6),\ (2,3,4,4,5,6),\ (2,3,3,4,6,6),\\
(2,3,3,4,5,7),\ (2,3,3,3,4,9),\ (2,2,5,5,5,5),\ (2,2,4,4,6,6),\ (2,2,4,4,4,8),\\
(2,2,3,5,5,7),\ (2,2,3,4,5,8),\ (2,2,3,3,7,7)^{\times2},\ (2,2,3,3,4,10),\ (2,2,2,6,6,6),\\
(2,2,2,4,6,8),\ (2,2,2,3,5,10),\ (2,2,2,3,3,12),\ (2,2,2,2,8,8),
\end{gathered}
\]
where $(2,2,3,3,7,7)^{\times2}$ indicates two combinatorially
distinct embeddings sharing that signature.
\end{theorem}

\begin{proof}[Proof (computational)]
By breadth-first search over \WM{} sequences, implemented at the dart
level with edge identities carried by the original cube labeling
$e_1,\dots,e_{12}$ propagated forward (a contracted edge's label is
inherited by the newly created edge); admissibility checks only
connectivity and absence of a monogon. Isomorphism of two results is
tested by a genuine orientation-preserving dart-system isomorphism
--- a bijection $\varphi$ conjugating $\alpha$ to $\alpha'$ and $\rho$
to $\rho'$ --- and \emph{not} by the coarser test of abstract graph
isomorphism, which discards $\rho$ and can wrongly identify
genuinely different embeddings; this distinction was confirmed by two
independent implementations of the isomorphism test (direct
propagation along $\alpha,\rho$, and an auxiliary directed
``flag-graph'' encoding checked by a general-purpose graph
isomorphism routine), which agreed on every instance tested. The
search terminates --- zero new isomorphism classes appear --- after
six \WM{}s from $d_1$.
\end{proof}

\begin{remark}[Independent confirmation]
This count of twenty was found independently, by an unrelated
enumeration method (contracting a perfect matching to a $4$-valent
medial graph and classifying those), in the author's
thesis~\cite[\S2.4]{Rodado2007}, where these are called the twenty
\emph{generic graphs} Gen1,\ldots,Gen20; we verified our own list
against the thesis's Table~2.1 face-degree by face-degree and found
exact agreement, including the coincidence that two combinatorially
distinct types share a signature. The thesis makes the same
embedded-versus-abstract distinction we make above: it separately
records that only seventeen of the twenty remain distinct as
\emph{abstract} (unembedded) graphs, since three pairs among the
twenty share an underlying abstract graph realized by two different
planar embeddings~\cite[\S7.1]{Rodado2007}.
\end{remark}

\begin{table}[h]
\centering
\small
\begin{tabular}{ccc}
\toprule
type & face-degree signature & shared? \\
\midrule
$g_0=d_1$ & $(4,4,4,4,4,4)$ & no \\
$g_1=d_2$ & $(3,3,4,4,5,5)$ & no \\
$g_8$ & $(2,2,3,3,7,7)$ & \textbf{yes} --- with $g_{19}$ \\
$g_{19}$ & $(2,2,3,3,7,7)$ & \textbf{yes} --- with $g_8$ \\
\bottomrule
\end{tabular}
\caption{Face-degree signature is a necessary but not sufficient
invariant of type: it distinguishes eighteen of the twenty types
outright, but $g_8$ and $g_{19}$ (both non-generic) share a signature
while being combinatorially distinct embedded graphs.}
\label{tab:collision}
\end{table}

\begin{remark}
Every one of the eighteen non-generic types has at least one bigon
face. By Lemma~\ref{lem:flippreservessimple} and its primal
counterpart, these are exactly the outcomes that Theorem~\ref{thm:orbit}
correctly excludes as inadmissible within the generic family: the
angle system, circle packing, and cell-complex structure of
Sections~\ref{sec:teich}--\ref{sec:circlepacking} are built on the
generic family alone, and Theorem~\ref{thm:fullorbit} records exactly
how much combinatorics is being set aside by that restriction.
\end{remark}

\begin{figure}[h]
\centering
\begin{tikzpicture}[vtx/.style={circle,draw,fill=white,minimum size=7mm,font=\scriptsize,inner sep=0pt},scale=0.85]
  \node[vtx] (c0) at (0.65,3.86) {$g_0$};
  \node[vtx] (c1) at (0.53,2.66) {$g_1$};
  \node[vtx] (c2) at (0.94,1.62) {$g_2$};
  \node[vtx] (c3) at (-0.55,1.86) {$g_3$};
  \node[vtx] (c4) at (-2.34,2.13) {$g_4$};
  \node[vtx] (c5) at (-0.72,0.68) {$g_5$};
  \node[vtx] (c6) at (0.10,-0.35) {$g_6$};
  \node[vtx] (c7) at (2.10,1.05) {$g_7$};
  \node[vtx] (c8) at (-2.30,0.90) {$g_8$};
  \node[vtx] (c9) at (-0.95,-0.65) {$g_9$};
  \node[vtx] (c10) at (-2.85,-0.35) {$g_{10}$};
  \node[vtx] (c11) at (-0.35,-1.55) {$g_{11}$};
  \node[vtx] (c12) at (2.35,-0.20) {$g_{12}$};
  \node[vtx] (c13) at (0.65,-2.35) {$g_{13}$};
  \node[vtx] (c14) at (-3.65,1.85) {$g_{14}$};
  \node[vtx] (c15) at (-1.35,-2.35) {$g_{15}$};
  \node[vtx] (c16) at (-0.55,-3.35) {$g_{16}$};
  \node[vtx] (c17) at (2.55,-1.55) {$g_{17}$};
  \node[vtx] (c18) at (-1.15,-4.35) {$g_{18}$};
  \node[vtx] (c19) at (3.75,-2.15) {$g_{19}$};
  \draw (c0)--(c1) node[midway,fill=white,inner sep=0.5pt,font=\tiny]{12};
  \draw (c1) to[loop right,looseness=8] node[right,font=\tiny]{4} (c1);
  \draw (c1)--(c2) node[midway,fill=white,inner sep=0.5pt,font=\tiny]{5};
  \draw (c1)--(c3) node[midway,fill=white,inner sep=0.5pt,font=\tiny]{3};
  \draw (c1)--(c4) node[midway,fill=white,inner sep=0.5pt,font=\tiny]{2};
  \draw (c2) to[loop above,looseness=8] node[above,font=\tiny]{2} (c2);
  \draw (c2)--(c3) node[midway,fill=white,inner sep=0.5pt,font=\tiny]{2};
  \draw (c2)--(c5) node[midway,fill=white,inner sep=0.5pt,font=\tiny]{3};
  \draw (c2)--(c6) node[midway,fill=white,inner sep=0.5pt,font=\tiny]{3};
  \draw (c2)--(c7) node[midway,fill=white,inner sep=0.5pt,font=\tiny]{4};
  \draw (c3)--(c5) node[midway,fill=white,inner sep=0.5pt,font=\tiny]{6};
  \draw (c3)--(c8) node[midway,fill=white,inner sep=0.5pt,font=\tiny]{3};
  \draw (c4)--(c6) node[midway,fill=white,inner sep=0.5pt,font=\tiny]{8};
  \draw (c4)--(c8) node[midway,fill=white,inner sep=0.5pt,font=\tiny]{3};
  \draw (c5)--(c9) node[midway,fill=white,inner sep=0.5pt,font=\tiny]{3};
  \draw (c5)--(c10) node[midway,fill=white,inner sep=0.5pt,font=\tiny]{6};
  \draw (c5)--(c11) node[midway,fill=white,inner sep=0.5pt,font=\tiny]{2};
  \draw (c6)--(c9) node[midway,fill=white,inner sep=0.5pt,font=\tiny]{2};
  \draw (c6)--(c11) node[midway,fill=white,inner sep=0.5pt,font=\tiny]{6};
  \draw (c6)--(c12) node[midway,fill=white,inner sep=0.5pt,font=\tiny]{5};
  \draw (c7)--(c12) node[midway,fill=white,inner sep=0.5pt,font=\tiny]{4};
  \draw (c8)--(c9) node[midway,fill=white,inner sep=0.5pt,font=\tiny]{6};
  \draw (c8)--(c14) node[midway,fill=white,inner sep=0.5pt,font=\tiny]{4};
  \draw (c9)--(c15) node[midway,fill=white,inner sep=0.5pt,font=\tiny]{1};
  \draw (c11)--(c16) node[midway,fill=white,inner sep=0.5pt,font=\tiny]{2};
  \draw (c12)--(c17) node[midway,fill=white,inner sep=0.5pt,font=\tiny]{2};
  \draw (c13) to[loop left,looseness=8] node[left,font=\tiny]{4} (c13);
  \draw (c13)--(c15) node[midway,fill=white,inner sep=0.5pt,font=\tiny]{2};
  \draw (c13)--(c16) node[midway,fill=white,inner sep=0.5pt,font=\tiny]{4};
  \draw (c13)--(c17) node[midway,fill=white,inner sep=0.5pt,font=\tiny]{4};
  \draw (c15) to[loop below,looseness=8] node[below,font=\tiny]{2} (c15);
  \draw (c15)--(c16) node[midway,fill=white,inner sep=0.5pt,font=\tiny]{4};
  \draw (c16) to[loop below,looseness=8] node[below,font=\tiny]{2} (c16);
  \draw (c16)--(c18) node[midway,fill=white,inner sep=0.5pt,font=\tiny]{2};
  \draw (c17) to[loop right,looseness=8] node[right,font=\tiny]{4} (c17);
  \draw (c17)--(c19) node[midway,fill=white,inner sep=0.5pt,font=\tiny]{2};
  \draw (c18) to[loop below,looseness=8] node[below,font=\tiny]{4} (c18);
  \draw (c19) to[loop right,looseness=8] node[right,font=\tiny]{6} (c19);
\end{tikzpicture}
\caption{The meta-graph of Theorem~\ref{thm:fullorbit}: vertices are
the twenty types, and an edge label $m$ between $g_i,g_j$ (or a
self-loop, for self-transitions) records that $m$ of the twelve
edges of the dual triangulation of $g_i$ admit a \WM{} landing on $g_j$;
vertex, including any self-loop, sum to at most twelve, the shortfall
being the count of inadmissible facets.}
\label{fig:metagraph}
\end{figure}
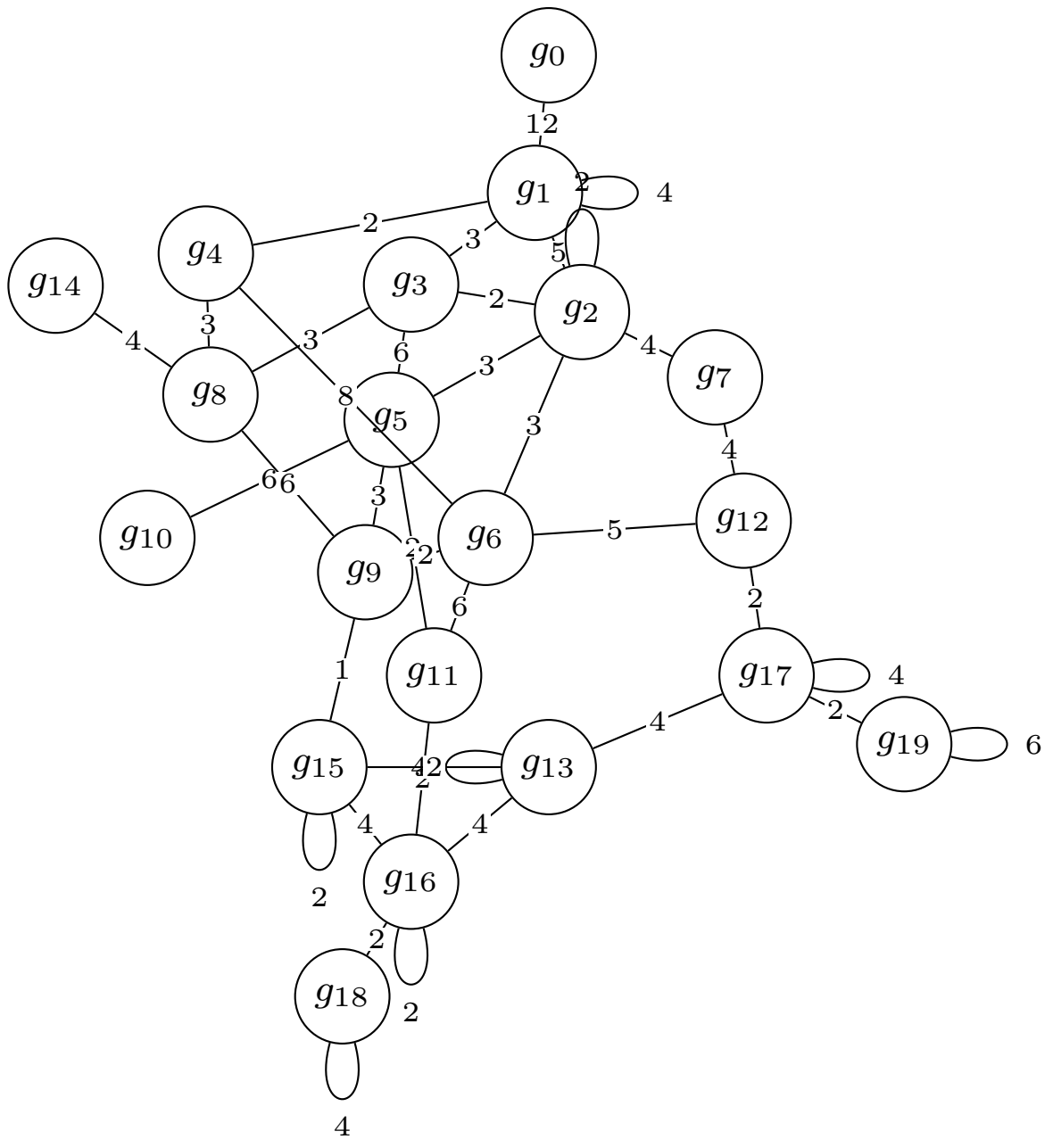

Figure~\ref{fig:small} shows the two generic types and their
Delaunay-dual triangulations explicitly; Figure~\ref{fig:metagraph}
displays the full twenty-type orbit as a meta-graph, whose vertices
are the isomorphism classes and whose edge multiplicities record how
many of a type's twelve facets realize a \WM{} to each neighboring
type (a self-loop recording self-transitions).

\begin{figure}[h]
\centering
\begin{tikzpicture}[scale=0.9,vtx/.style={circle,fill=black,inner sep=1.6pt},
  every node/.style={font=\scriptsize}]
  \coordinate (Fppp) at (2,2);
  \coordinate (Fmpp) at (-2,2);
  \coordinate (Fmmp) at (-2,-2);
  \coordinate (Fpmp) at (2,-2);
  \coordinate (Fppm) at (1,1);
  \coordinate (Fmpm) at (-1,1);
  \coordinate (Fmmm) at (-1,-1);
  \coordinate (Fpmm) at (1,-1);
  \draw (Fppp)--(Fmpp)--(Fmmp)--(Fpmp)--cycle;
  \draw (Fppm)--(Fmpm)--(Fmmm)--(Fpmm)--cycle;
  \draw (Fppp)--(Fppm); \draw (Fmpp)--(Fmpm);
  \draw (Fmmp)--(Fmmm); \draw (Fpmp)--(Fpmm);
  \foreach \c in {Fppp,Fmpp,Fmmp,Fpmp,Fppm,Fmpm,Fmmm,Fpmm}
    {\node[vtx] at (\c) {};}
\end{tikzpicture}
\hspace{1cm}
\begin{tikzpicture}[scale=0.9,vtx/.style={circle,fill=black,inner sep=1.6pt},
  every node/.style={font=\scriptsize}]
  \coordinate (v1) at (0,2);
  \coordinate (v2) at (-2,-1.2);
  \coordinate (v3) at (2,-1.2);
  \coordinate (v1p) at (0,-0.6);
  \coordinate (v2p) at (0.8,0.3);
  \coordinate (v3p) at (-0.8,0.3);
  \draw (v1)--(v2)--(v3)--cycle;
  \draw (v1p)--(v2p)--(v3p)--cycle;
  \draw (v1)--(v2p); \draw (v1)--(v3p);
  \draw (v2)--(v1p); \draw (v2)--(v3p);
  \draw (v3)--(v1p); \draw (v3)--(v2p);
  \foreach \c/\lab in {v1/$1$,v2/$2$,v3/$3$,v1p/$1'$,v2p/$2'$,v3p/$3'$}
    {\node[vtx,label=\lab] at (\c) {};}
\end{tikzpicture}
\caption{Left: the cube graph $d_1$ (a Schlegel diagram). Right: its
dual, the octahedral triangulation $T$, on the six marked points
with antipodal (non-adjacent) pairs $\{1,1'\},\{2,2'\},\{3,3'\}$.}
\label{fig:small}
\end{figure}
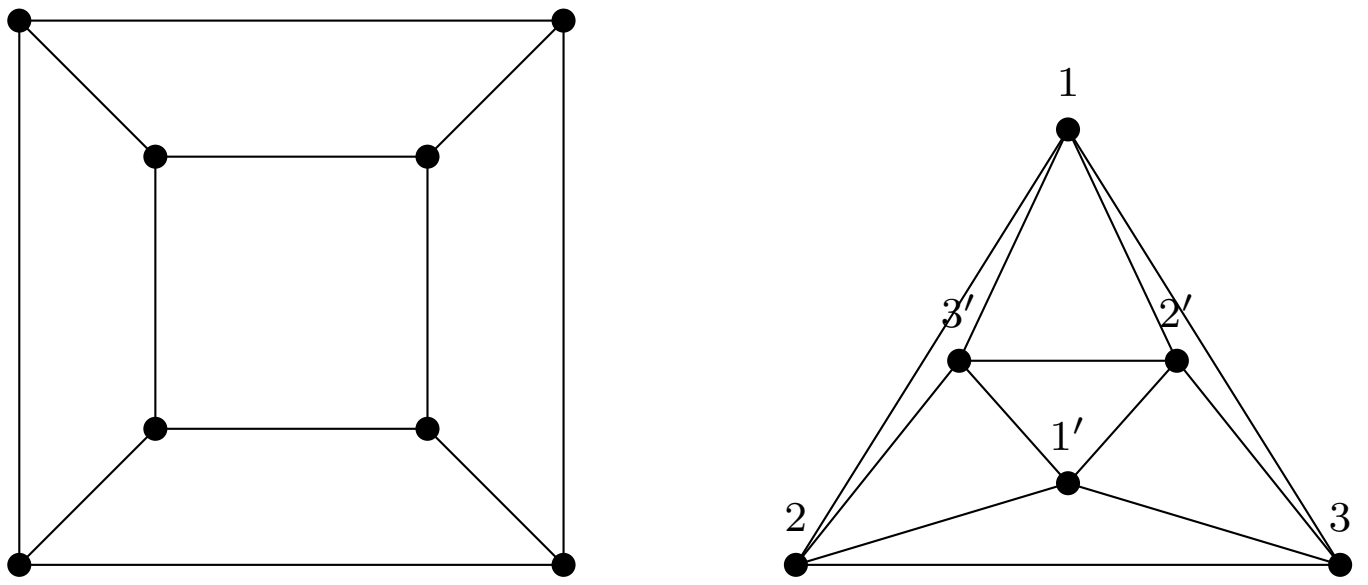

\subsection{The contraction hierarchy}
\label{sec:hierarchy}

The bigons of Section~\ref{sec:fullorbit} leave the vertex count fixed
at eight: a \WM{} contracts an edge and immediately re-expands it
along the opposite diagonal, so $V=8$, $E=12$ throughout. A different,
and in fact more fundamental, generalization is obtained by
\emph{not} re-expanding: contracting a single edge of a generic
CE-graph merges its two trivalent endpoints into one $4$-valent
vertex, leaving a graph with seven vertices, eleven edges, and ---
automatically, by Euler's formula, provided the result stays
connected and monogon-free --- still six faces. Iterating drives the
vertex count down one at a time, and reveals the complete face
structure of the local cell complex that the generic trivalent level
merely sits at the top of.

\begin{definition}[Contraction hierarchy]
\label{def:hierarchy}
For a CE-type graph $G$ (connected, cellularly embedded, six faces,
no monogon --- trivalence is \emph{not} required beyond level $0$)
and an edge $e=(p,q)$ of $G$, the \emph{contraction} $G/e$ merges $p,q$
into a single vertex of degree $\deg(p)+\deg(q)-2$, deleting $e$. The
contraction is \emph{admissible} exactly when $G/e$ has no monogon.
Level $0$ is the generic trivalent family of Sections~\ref{sec:classification}--\ref{sec:orbit};
level $k+1$ is the set of graphs obtainable from a level-$k$ graph by
one admissible contraction.
\end{definition}

\begin{lemma}
\label{lem:levelinvariant}
Every graph at level $k$ of the hierarchy has $8-k$ vertices,
$12-k$ edges, and (still) exactly six faces.
\end{lemma}

\begin{proof}
Contraction always removes exactly one vertex and one edge; the claim
for $k=0$ is the base case (Section~\ref{sec:classification}), and
$V-E+F=2$ together with $\Delta V=\Delta E=-1$ forces $\Delta F=0$ at
every admissible step.
\end{proof}

A subtlety not visible at level $0$, where every graph is simple,
must be handled once vertices of degree $\ge 4$ appear: a level-$k$
graph can carry a self-loop (an edge both of whose darts sit at the
same vertex) without that loop bounding a monogon, if enough further
structure survives on the loop's far side. Contracting such a loop is
not the operation of Definition~\ref{def:hierarchy} at all --- it
would delete an edge without merging two vertices, changing $E$ but
not $V$ --- so loops are simply excluded from the set of edges a
level-$k$ graph offers for the next contraction; every non-loop edge
of a monogon-free graph remains eligible.

\begin{theorem}[The complete contraction hierarchy]
\label{thm:hierarchy}
Up to embedded graph equivalence, the contraction hierarchy of
Definition~\ref{def:hierarchy}, started from the two generic types,
has exactly seven levels, $k=0,\dots,6$, with
\[
\begin{array}{c|ccccccc}
k & 0 & 1 & 2 & 3 & 4 & 5 & 6\\
\hline
(V,E) & (8,12) & (7,11) & (6,10) & (5,9) & (4,8) & (3,7) & (2,6)\\
\hline
\#\text{types} & 2 & 6 & 22 & 32 & 22 & 5 & 1
\end{array}
\]
(ninety types in total), and level $6$ is a single graph
$\Theta_6$: two vertices joined by six parallel edges, with all six
faces bigons.
\end{theorem}

\begin{proof}[Proof (computational)]
By level-by-level breadth-first search: at each level, every non-loop
edge of every representative graph is contracted (via the dart-level
splice of Definition~\ref{def:hierarchy}, generalized from the
trivalent case by allowing $\deg(p),\deg(q)\ge 3$ arbitrary), the
monogon-free results are kept, and duplicates are removed by the same
genuine orientation-preserving dart-system isomorphism used
throughout this paper (never by the coarser, rotation-forgetting
abstract graph isomorphism). Every contraction is checked to lower
$V$ by exactly one, which is why loop edges must be excluded. The search
terminates at $V=2$ with a single isomorphism class, matching
Lemma~\ref{lem:levelinvariant}'s prediction $(V,E)=(2,6)$ exactly.
\end{proof}

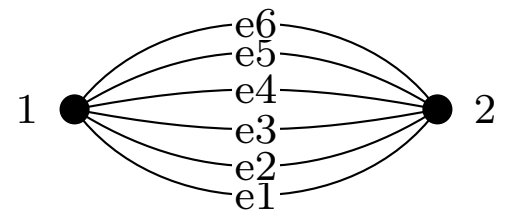
\begin{figure}[h]
\centering
\begin{tikzpicture}[vtx/.style={circle,fill=black,inner sep=2pt},every node/.style={font=\scriptsize}]
  \coordinate (L) at (-1.2,0);
  \coordinate (R) at (1.2,0);
  \foreach \rad/\lab in {-0.9/e1,-0.54/e2,-0.18/e3,0.18/e4,0.54/e5,0.9/e6}
    {\draw (L) to[bend left=\rad*60] node[midway,sloped,fill=white,inner sep=0.5pt]{\lab} (R);}
  \node[vtx,label=left:$1$] at (L) {};
  \node[vtx,label=right:$2$] at (R) {};
\end{tikzpicture}
\caption{The unique deepest-level graph $\Theta_6$: two vertices of
degree six, joined by six parallel edges, bounding six bigon faces.
Every one of the ninety types across all seven levels of
Theorem~\ref{thm:hierarchy} eventually contracts down to this single
graph.}
\label{fig:theta6}
\end{figure}

\begin{remark}[Geometric meaning of a contraction]
\label{rem:contraction-meaning}
Dually, contracting an edge of $G$ deletes the corresponding edge of
$T$, merging two adjacent Delaunay triangles into a single
quadrilateral. This is exactly the classical Delaunay wall-crossing
degeneracy: the shared diagonal becomes removable precisely when its
four endpoints become concyclic (Proposition~\ref{prop:xezero}). A
level-$k$ graph therefore records $k$ simultaneous concyclicity
coincidences among the six marked points. They are independent, so that
the stratum has codimension $k$, unless the dual triangulation $T$ is
bipartite (every vertex of $G$ has even degree): then one condition
follows from the others (Lemma~\ref{lem:incidencerank}) and the
codimension is $k-1$. The contraction hierarchy is a combinatorial
model, in the spirit of the face lattice of a secondary
polytope~\cite{GKZ1994}, of the strata where the generic (simplicial)
Delaunay decomposition degenerates; only the realizable types occur as
strata (Remark~\ref{rem:combnotreal}). The deepest
stratum, $\Theta_6$, is the fully degenerate configuration in which
every Delaunay cell has merged into one: a single circle through all
six points remains as the only trace of the original triangulation,
subdivided into six arcs by the six marked points themselves, dual to
the six parallel edges of Figure~\ref{fig:theta6}.
\end{remark}

\begin{remark}[$\Theta_6$ is realized]
\label{rmk:theta6realized}
Unlike the general hierarchy (Remark~\ref{rem:combnotreal}), this
deepest stratum is not merely combinatorially conceivable: it is
geometrically realized, and admits an essentially closed-form
construction bypassing the degenerate-limit numerics of
Section~\ref{sec:contractrealize} entirely. Since $T$ has only six
vertices, six edges, and (dually) two faces, Euler's formula forces
$T$ itself to be a bare hexagonal cycle, so the whole sphere is two
copies of a single hyperbolic hexagon glued along their common
boundary; for the cone angle to be exactly $\pi$ at each of the six
marked points (half contributed by each copy), the hexagon must be
\emph{right-angled}, a classical object with an explicit closed-form
relation between alternating side lengths~(see e.g.~\cite{Buser1992}).
For the six points to be concyclic, as $\Theta_6$ requires, the hexagon
must moreover be inscribed in a circle. Right-angled hexagons form a
three-parameter family and the inscribed ones a one-parameter family:
the cell of $\Theta_6$ is one-dimensional, as Lemma~\ref{lem:incidencerank}
predicts for the bipartite hexagon $T$. For one such hexagon we verified
directly (closure of the hexagon in coordinates, all six
angles equal to $\pi/2$, a common circumcenter equidistant from all
six vertices confirming the two hexagon-centers are genuine
degree-six Voronoi vertices matching $\Theta_6$'s own vertex degrees,
and Gauss--Bonnet giving area $2\pi$) that this construction is
completely consistent. The thesis reaches the same conclusion by
a related but distinct explicit construction --- the surface itself
(not just its orbifold quotient) built directly from \emph{four}
right-angled hyperbolic hexagons glued along a dodecagon, rather than
two hexagons glued to form the quotient which is then doubled by the
branched cover as in Theorem~\ref{thm:realizability} --- and
identifies a one-real-parameter family of such
surfaces~\cite[\S2.5.1]{Rodado2007}; the two constructions agree on
the total area ($4\pi$ for the closed genus-two surface either way)
and on the combinatorics, and we regard them as independent
confirmations of the same fact.
\end{remark}

\subsection{Level-by-level atlas of the hierarchy}
\label{sec:atlas}

\begin{figure}[p]
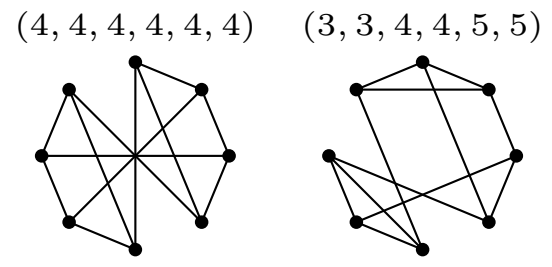

\centering

\caption{Level 0 of the contraction hierarchy: $V=8$, $E=12$, 2 types, up to embedded graph equivalence, each labeled by its face-degree signature.}
\label{fig:level0}
\end{figure}

Figures~\ref{fig:level0}--\ref{fig:level6} display all ninety types
explicitly, one table per level, each graph drawn on its own
canonical circular vertex layout (Section~\ref{sec:hierarchy}'s
proof) and labeled by its face-degree signature; a bundle of
parallel arcs between two vertices records a multiple edge of that
multiplicity (up to six, at $\Theta_6$ itself). No two graphs drawn
within a single level are isomorphic; graphs across different levels
are of course never isomorphic either, since they have different
vertex counts.

\begin{figure}[p]
\centering
%
\caption{Level 1 of the contraction hierarchy: $V=7$, $E=11$, 6 types, up to embedded graph equivalence, each labeled by its face-degree signature.}
\label{fig:level1}
\end{figure}

\begin{figure}[p]
\centering
%
\caption{Level 2 of the contraction hierarchy: $V=6$, $E=10$, 22 types, up to embedded graph equivalence, each labeled by its face-degree signature.}
\label{fig:level2}
\end{figure}

\begin{figure}[p]
\centering
%
\caption{Level 6 of the contraction hierarchy: $V=2$, $E=6$, 1 type, up to embedded graph equivalence, each labeled by its face-degree signature.}
\label{fig:level6}
\end{figure}

\subsection{Regenerating the hierarchy from the deepest graph}
\label{sec:regeneration}

The contraction hierarchy was built downward, from the generic level
to $\Theta_6$. We now ask whether it can be built back \emph{upward},
starting only from $\Theta_6$, and if so, by what moves.

\begin{definition}[Expansion]
\label{def:expansion}
The \emph{expansion} of a graph $G$ at a vertex $v$ of degree $d\ge4$,
along a contiguous arc of $2\le a\le d-2$ darts in the cyclic rotation
at $v$, is the precise reverse of Definition~\ref{def:hierarchy}'s
contraction: the arc's darts and the complementary $(d-a)$ darts
become two new vertices of degrees $a+1$ and $d-a+1$, joined by one
new edge that closes each of their rotations. Contracting the new
edge recovers $G$ exactly.
\end{definition}

Unlike contraction, expansion is highly non-unique: a degree-$d$
vertex admits, up to the trivial symmetry of which side is called
which, on the order of $d^2$ inequivalent arcs to expand along. We
verified Definition~\ref{def:expansion} computationally by checking
both round trips on many examples: contracting an edge and then
expanding the resulting vertex along some arc always recovers the
original graph, and expanding a vertex and then contracting the new
edge always recovers the pre-expansion graph exactly, not merely up
to isomorphism.

\begin{theorem}[Regeneration]
\label{thm:regeneration}
Starting from $\Theta_6$ alone and closing under \emph{both}
expansion and contraction moves (discarding any result with more
than eight vertices or with a monogon), the reachable family
converges, after eleven rounds, to exactly $280$ types. It contains
all ninety types of Theorem~\ref{thm:hierarchy}, and its restriction
to eight vertices is exactly the twenty types of
Theorem~\ref{thm:fullorbit} -- an independent, purely bottom-up
confirmation of that top-down result.
\end{theorem}

\begin{proof}[Proof (computational)]
Breadth-first search from $\Theta_6$, alternately applying every
admissible contraction (Definition~\ref{def:hierarchy}) and every
expansion (Definition~\ref{def:expansion}, over every vertex of
degree $\ge4$ and every valid arc) to every representative found so
far, deduplicated at each round by the same orientation-preserving
dart-system isomorphism used throughout. The count stabilizes at
$280$ (zero new classes in round twelve). Matching each of the
ninety generic-rooted types of Theorem~\ref{thm:hierarchy} against
this pool of $280$ finds all ninety present; separately, filtering
the pool to its $20$ eight-vertex members and matching against
Theorem~\ref{thm:fullorbit}'s twenty types finds all twenty present,
with no extra eight-vertex types beyond them.
\end{proof}

The full $280$-type closure is larger than the $90$-type generic-rooted
hierarchy because it also contains the eighteen non-generic
trivalent types of Theorem~\ref{thm:fullorbit} together with
\emph{their own} contraction sub-hierarchies, which Theorem~\ref{thm:hierarchy}
never had occasion to compute, having started only from the two
generic types. Theorem~\ref{thm:regeneration} shows this is exactly
the same object either way: whatever is reachable downward from the
generic level by contraction is also reachable upward from $\Theta_6$
by a mix of expansion and contraction.

\begin{proposition}[Pure expansion does not suffice]
\label{prop:noexpansiononly}
Closing under expansion \emph{alone}, starting from $\Theta_6$
(discarding results with more than eight vertices), converges after
six rounds to exactly $128$ types. Both generic types $d_1,d_2$ are
among them, but only sixty-six of Theorem~\ref{thm:hierarchy}'s
ninety generic-rooted types are reached; the other twenty-four ---
none of them at level $0$ or level $6$, all at intermediate levels
--- are never produced by expansion alone, no matter how it is
iterated.
\end{proposition}

\begin{proof}[Proof (computational)]
Same breadth-first search as Theorem~\ref{thm:regeneration}, with the
contraction moves omitted. The closure stabilizes at $128$ (round
seven adds nothing). Matching against the ninety generic-rooted types
finds sixty-six present and twenty-four absent; $d_1$ and $d_2$ are
directly verified to be among the sixty-six.
\end{proof}

Proposition~\ref{prop:noexpansiononly} shows that $\Theta_6$ is the
unique sink of the contraction hierarchy (every type contracts to it, by
Theorem~\ref{thm:hierarchy}) but not a source under expansion alone:
some intermediate types are reached only along paths that expand,
contract and expand again. The whole family is connected only when
moves in both directions are allowed.

\begin{remark}[Independent verification]
\label{rem:indep}
All enumerations in this section were recomputed by a second,
independent program that works with the dual triangulation $T$ rather
than with dart systems of $G$. A contraction of a non-loop edge of $G$
is the deletion of a non-bridge edge of $T$, an expansion is the
insertion of a diagonal in a face of $T$ of degree at least four, and
the no-monogon condition says that $T$ has no vertex of degree one.
Maps are compared by a canonical breadth-first code of the rotation
system. This reproduces the twenty trivalent types of
Theorem~\ref{thm:fullorbit}, with exactly seventeen distinct abstract
graphs as in~\cite[\S7.1]{Rodado2007}; the counts $2,6,22,32,22,5,1$ of
Theorem~\ref{thm:hierarchy}; the closure of $280$ types of
Theorem~\ref{thm:regeneration}; and the $128$ types, with $66$ of the
ninety reached and $24$ missed, of
Proposition~\ref{prop:noexpansiononly}. The missed types lie at levels
$2,3,4,5$, numbering $2,9,10,3$. The program also identifies the
closure of Theorem~\ref{thm:regeneration}: its $280$ types are exactly
all maps on the sphere with six vertices and no vertex of degree one,
of which $20,55,97,71,31,5,1$ have $12,11,\dots,6$ edges.
Testing each type for a strictly positive solution of the vertex
equations by linear programming (the linear part of the criterion of
Section~\ref{sec:springborncorrect}) gives $10$ of the $20$ trivalent
types, as in~\cite[Ch.~4]{Rodado2007}; $116$ of the $280$ types,
including $26$ with eleven edges and $2$ with seven, as
in~\cite[Props.~7.2 and~7.6]{Rodado2007}; and $2,5,13,12,6,2,1$ types
by level in the ninety-type hierarchy, $41$ in all.
\end{remark}

\section{Discussion and further directions}
\label{sec:discussion}

\subsection{Relation to decorated Teichm\"uller theory}

The Bowditch--Epstein/Penner construction~\cite{Penner1987,BowditchEpstein1988}
attaches a canonical ideal-triangulation cell decomposition to
\emph{any} decorated hyperbolic surface with cusps, with cells glued
along diagonal flips exactly as in Section~\ref{sec:orbit}. The
present construction is the closed, cone-point analogue: rather than
horocyclic decorations at cusps, we use the intrinsic cone angle
$\pi$ forced at each of the six branch points by
Proposition~\ref{prop:coneangle}. Theorem~\ref{thm:classification}
and Theorem~\ref{thm:orbit} show that in genus two this general
machinery collapses to an unusually small and completely explicit
structure --- two cells and a five-edge transition graph on the
stacked cell --- in sharp contrast to the generally infinite flip
graphs that arise for higher-complexity punctured surfaces. It is
natural to ask whether comparably small and complete flip-graph
classifications occur for other low-complexity closed or once-marked
hyperbolic orbifolds with cone points of order $2$, and if so whether
their size grows in a controlled way with the orbifold Euler
characteristic. In this dictionary, Section~\ref{sec:hierarchy}'s
contraction hierarchy is the more fundamental object: in the
Bowditch--Epstein/Penner cell complex, the top-dimensional cells are
indexed by ideal triangulations (our trivalent level $0$), and each
codimension-$k$ face of the complex is indexed by a partial
degeneration in which $k$ edges have been contracted -- exactly
Definition~\ref{def:hierarchy}. Theorem~\ref{thm:hierarchy} therefore
describes the combinatorial face structure below the two generic types;
the cells that actually occur are the realizable ones, $41$ of the
ninety (Remark~\ref{rem:indep}). The complete decomposition, obtained by
contracting from all ten realizable trivalent types, is computed
in~\cite[Ch.~2,~6--7]{Rodado2007}, which also determines realizability
at every level; it has $116$ cells, among the $280$ types of
Theorem~\ref{thm:regeneration}. The finer question of whether the pattern
$2,6,22,32,22,5,1$ -- in particular its symmetry about level $3$,
$22$ types appearing at both level $2$ and level $4$ -- reflects a
genuine duality of this narrower hierarchy (analogous to the
self-duality of the face lattice of a simple polytope) or is
coincidental to genus two, and to this specific two-generator
sub-hierarchy, remains open.

\subsection{Relation to circle patterns and discrete uniformization}

Theorem~\ref{thm:rankAi} establishes agreement of dimension between
the linear cone-angle system and $\TT_2$, but the linear system of
Proposition~\ref{prop:linearconstraint} is a relaxation of the true,
generally nonlinear circle-packing equations governed by the
hyperbolic law of cosines~\cite{Rivin1994,ThurstonShapes1998}.
Promoting Theorem~\ref{thm:rankAi} to an actual local diffeomorphism
$D_i^\circ \to \TT_2$ requires the existence and local rigidity theory
for circle patterns with prescribed combinatorics and intersection
angles on a fixed conformal class. Section~\ref{sec:realizability}
carries out the existence half of this explicitly for both generic
types, directly from the hyperbolic law of cosines rather than any
named existence theorem: we solve for edge lengths realizing each
$T_i$ as an honest hyperbolic triangulation with the correct cone
angle at every vertex (exactly, for $d_1$, by a symmetric closed form;
numerically, to machine precision, for $d_2$), verify the result is
genuinely the Delaunay decomposition via the local empty-circumcircle
criterion, and recover Proposition~\ref{prop:linearconstraint}'s
$x_e$ from the resulting circumcenters. This produces one explicit
hyperbolic sphere per generic type, and hence (via
Section~\ref{sec:realizability}'s cutting-and-doubling construction)
one actual genus-two hyperbolic surface realizing each of $d_1,d_2$.
We further extend this one step into the contraction hierarchy of
Section~\ref{sec:hierarchy} (Theorem~\ref{thm:contractrealize}),
verifying directly for one edge of each generic type that the linear
system's boundary stratum $x_e=0$ is reached by an honest limit of
this same family, and correctly identifying what that limit means
geometrically (Proposition~\ref{prop:xezero}) along the way. Whether
this can be extended from a single point to the whole of
$D_i^\circ$, promoted to a local diffeomorphism (injective with
nonsingular derivative) rather than a single realized point per type,
and carried systematically through the rest of the ninety-type
hierarchy, remains for future work.

\subsection{Relation to the equisymmetric stratification of $\MM_2$}

Automorphisms of a genus-two curve $X$ correspond, via
Lemma~\ref{lem:allhyp}, to M\"obius transformations of $\PP^1$
permuting the six branch points; the classical classification of
genus-two curves with extra automorphisms goes back to
Bolza~\cite{Bolza1887} and has been organized modernly into the
equisymmetric stratification of $\MM_2$ (see e.g.\
\cite{Broughton1991}). The maximally symmetric genus-two curve, the
Bolza curve $y^2=x^5-x$, has automorphism group of order $48$; its six
branch points $0,\infty,\pm1,\pm i$ are the vertices of a regular
octahedron inscribed in $\PP^1$. This is suggestively close to our type
$d_1$, whose orientation-preserving automorphism group is exactly the
order-$24$ rotation group of the octahedron
(Proposition~\ref{prop:aut}). A generic point of $\MM_2$, whose only
automorphism is the hyperelliptic involution, has a Delaunay
decomposition of type $d_2$; the order-two group $\Aut^+(d_2)$ is a
symmetry of the combinatorial type only, since the hyperelliptic
involution acts trivially on the quotient sphere and hence on $T$. The
expectation that the locus of type $d_1$ contains the Bolza curve is
settled in~\cite[\S8.3.1]{Rodado2007}, which identifies the Bolza curve
by its Delaunay angle coordinates as the point of $D_1$ with all twelve
interior angles equal to $3\pi/4$. This is the unique point of $D_1$
fixed by the automorphism group of $T_1$, which acts transitively on
the edges. In the exterior-angle normalization used here it is
$x_e=1-\tfrac34=\tfrac14$ for every edge, exactly the symmetric point of
Proposition~\ref{prop:cubeconstruction}, obtained in this paper from the
equilateral triangle equation alone. Relating the combinatorial
automorphism group $\Aut^+(d_i)$ to the geometric automorphism group of
the corresponding points of $\MM_2$ in general, beyond this case,
remains an appealing question for future work.

\subsection{Summary of what remains open}

We highlight two resolved steps, one partially resolved step, and two
remaining open steps suggested by this work:
\begin{enumerate}[label=(\roman*)]
\item \emph{Resolved in Section~\ref{sec:realizability}.} We construct,
directly from the hyperbolic law of cosines and verified independently
via a Gauss--Bonnet cross-check and a local Delaunay criterion, an
explicit hyperbolic sphere with six cone points realizing each generic
type's Delaunay decomposition (exactly, for $d_1$; numerically, to
machine precision, for $d_2$); this upgrades Theorem~\ref{thm:rankAi}
to an actual realization: both $d_1,d_2$ occur as the Delaunay/Voronoi
graph of a genuine genus-two hyperbolic surface. The thesis confirms
this, and much more, for the generic level as a whole: applying
Springborn's theorem with $\Sigma=T$
(Section~\ref{sec:springborncorrect}) and a linear-programming closing
condition, exactly ten of the twenty types of
Theorem~\ref{thm:fullorbit} are Delaunay
realizable~\cite[Ch.~4--5]{Rodado2007}, of which $d_1,d_2$ are two,
confirmed here by explicit construction;
\item \emph{Resolved by citation.} The natural expectation that
$d_1$'s locus meets the equisymmetric stratum of the Bolza curve is
already established in~\cite[\S8.3.1]{Rodado2007}, and matches, to
the exact numerical value, our own independently derived symmetric
solution above;
\item \emph{Partially resolved in Section~\ref{sec:contractrealize}.}
The boundary $x_e=0$ of the linear system is the concyclicity of the
four vertices of the two triangles at $e$
(Proposition~\ref{prop:xezero}), and for one admissible edge of each
generic type this limit is reached by an explicit hyperbolic
construction. Which contractions are realizable at all is answered
completely in~\cite[Ch.~7]{Rodado2007} (Remark~\ref{rem:combnotreal}).
What remains open is whether the explicit construction can be carried
out at every realizable instance, and iterated to $\Theta_6$ along every
realizable path, not just the closed-form one of
Remark~\ref{rmk:theta6realized};
\item determine precisely the relationship, discussed informally
above, between $\Aut^+(d_i)$ and the equisymmetric stratification of
$\MM_2$ beyond the single confirmed case of $d_1$ and the Bolza curve;
\item determine whether the eighteen non-generic types of
Theorem~\ref{thm:fullorbit}, and the deeper contraction hierarchy of
Section~\ref{sec:hierarchy}, correspond to walls of a compactified or
extended cell decomposition of $\TT_2$, in analogy with the way
degenerating Delaunay decompositions signal approach to strata of
$\overline{\MM}_2$ in the decorated-Teichm\"uller-theoretic
literature; item~(iii)'s verified instances, and the complete
realizability determination, are concrete data toward this picture.
\end{enumerate}

\section{Conclusion}

We have classified the Delaunay decompositions of the hyperelliptic
quotient of a genus-two Riemann surface: two generic types, connected by
a fully determined \WM{} orbit, whose linear angle system has rank
matching the six real dimensions of $\TT_2$. Both generic types are
realized by explicit hyperbolic cone spheres, obtained directly from the
hyperbolic law of cosines and checked by Gauss--Bonnet and the local
Delaunay criterion, and the corresponding genus-two surfaces are built
by cutting along a perfect matching of geodesics and doubling. We
recorded the more efficient route through Springborn's theorem with the
cell decomposition $\Sigma=T$, as in the author's
thesis~\cite{Rodado2007}, whose linear-programming treatment shows that
ten of the twenty trivalent types are realizable. The boundary $x_e=0$
of the angle system is the concyclicity of four branch points, and for
one edge of each generic type the corresponding contracted graph is
realized by an explicit limit.

Below the generic level, allowing bigons gives exactly twenty trivalent
types, agreeing signature for signature with the thesis, and the
contraction hierarchy from the two generic types has
$2,6,22,32,22,5,1$ types, ninety in all, ending at $\Theta_6$ (the graph
$I(6)$ of the thesis), for which we give a closed-form realization. The
hierarchy can be regenerated from $\Theta_6$ using both expansion and
contraction, which yields all $280$ maps on the sphere with six vertices
and no vertex of degree one, but not by expansion alone. The finite
combinatorics is formally verified in Lean, and every enumeration is
reproduced by an independent program. Genus two is distinguished among
closed surfaces by hyperellipticity, which is what makes this explicit,
finite control possible; we hope the methods and the open steps of
Section~\ref{sec:discussion} are useful for other low-complexity moduli
problems.

\section*{Code availability}

The Lean~4 formalization of Remark~\ref{rem:lean} and the enumeration
programs of Remark~\ref{rem:indep} are available at
\url{https://github.com/armandorodado/weierstrass-delaunay} and archived at
Zenodo, \href{https://doi.org/10.5281/zenodo.23870958}{doi:10.5281/10.5281/zenodo.23870958}

\end{document}